 \documentclass[final]{siamltex}

\usepackage{pst-node}
\usepackage{booktabs}
\usepackage{lscape}
\usepackage{array}
\usepackage{boxedminipage}
\usepackage{longtable}
\usepackage{graphicx}
\usepackage{psfrag}
\usepackage[small]{subfigure}
\usepackage{calc}

\usepackage{dcolumn}
\newcolumntype{d}[1]{D{.}{.}{#1}}

\input{macros.tex}
\renewcommand{\min}{\hbox{\normalfont  min}}
\renewcommand{\max}{\hbox{\normalfont  max}}
\newcommand{\aref}[1]{\textup{\ref{#1}}}
\newcommand{\mcol}[3]{\multicolumn{#1}{#2}{#3}}

\newcommand{\AugL}{\Lscr}

\newcommand{\dystark}{\Deltait\ystark}
\newcommand{\fomult}{\yhat}
\newcommand{\lsmult}{\ytilde}
\newcommand{\clink}[1]{\overline{c}\k(#1)}

\newcommand{\assign}{\leftarrow}
\newcommand{\iset}[2]{#1,\ldots,#2}

\newcommand{\pospart}[1]{[#1]^+}
\newcommand{\False}{\textsf{false}}
\newcommand{\True}{\textsf{true}}
\newcommand{\NP}  {{\rm(NP)}}
\newcommand{\NPi} {{\rm(NPi)}}
\newcommand{\GNP} {{\rm(GNP)}}
\newcommand{\ELCk}  {{\rm(ELC$\k$)}}

\newcommand{\ELCkp} {{\rm(ELC$\k'$)}}

\newcommand{\ELCik}{{\rm(ELCi$\k$)}}
\newcommand{\LCk} {{\rm(LC$\k$)}}

\newcommand{\BCk} {{\rm(BC$\k$)}}

\newcommand{\PP}[1]{{\rm(PP#1)}}

\newcommand{\seq}[1]{\ensuremath{\{#1\}}}
\newcommand{\sigmaB}{\overline{\sigma}}
\newcommand{\sigmab}{\underline{\sigma}}
\newcommand{\prob}[1]{{\it #1}}
\newcommand{\opt}[1]{{\tt #1}}
\newcommand{\LCLOPT}{LCLOPT}

\newtheorem{assumption}[theorem]{Assumption}

\usepackage[algonl,vlined,ruled]{algorithmMPF}
\newenvironment{algo}[1]
{\def\theAlgoLine{\arabic{AlgoLine}}
 \begin{algorithm}[#1]%
    \SetCommBlock{Comment}
    \SetArgSty{texttsf}%
    \SetTitleSty{textsf}{}%
    \SetKwInput{Inputs}{Input}%
    \SetKwInput{Outputs}{Output}%
    \SetKwData{Converged}{converged}%
  }
  {\end{algorithm}}

\title{A globally convergent LCL method
    \\ for nonlinear optimization}

\author{%
  Michael P. Friedlander%
  \thanks{Mathematics and Computer Science Division,
          Argonne National Laboratory,
          9700 S. Cass Avenue,
          Argonne, IL 60439-4844 ({\tt michael@mcs.anl.gov}).
          This work was supported by the U.S. National Science
          Foundation grant CCR-9988205
          and by
          the Mathematical, Information, and Computational Sciences
          Division subprogram of the Office of Advanced Scientific
          Computing Research, U.S. Department of Energy
          contract W-31-109-Eng-38.}
  \and%
  Michael A. Saunders%
  \thanks{Department of Management Science and Engineering,
          Stanford University,
          Stanford, CA 94305-4026 ({\tt saunders@stanford.edu}).
          This work was supported by the U.S. National Science
          Foundation grant CCR-9988205, and the U.S. Office of Naval
          Research grants N00014-96-1-0274 and N00014-02-1-0076.
          \hfill Version~of~\hbox{\today}}}

\begin{document}

\newsavebox{\fminibox}
\begin{lrbox}{\fminibox}
\begin{minipage}{\linewidth-2\fboxsep-2\fboxrule}
{\bf Preprint ANL/MCS-P1015-1202, December, 2002}
\\   Mathematics and Computer Science Division
\\   Argonne National Laboratory
\end{minipage}\end{lrbox}

\vspace*{-3\baselineskip}\noindent\fbox{\usebox{\fminibox}}\vspace*{\baselineskip}

\maketitle

\begin{abstract}
  For optimization problems with nonlinear constraints, linearly
  constrained Lagrangian (LCL) methods sequentially minimize a
  Lagrangian function subject to linearized constraints.  These
  methods converge rapidly near a solution but may not be reliable
  from arbitrary starting points.  The well known example \MINOS\ has
  proven effective on many large problems.  Its success motivates us
  to propose a globally convergent variant.  Our stabilized LCL method
  possesses two important properties: the subproblems are always
  feasible, and they may be solved inexactly. These features are
  present in \MINOS\ only as heuristics.
  
  The new algorithm has been implemented in \Matlab, with the option
  to use either the \MINOS\ or \SNOPT\ Fortran codes to solve the
  linearly constrained subproblems.  Only first derivatives are
  required.  We present numerical results on a nonlinear subset of the
  \COPS, \CUTE, and HS test-problem sets, which include many large
  examples.  The results demonstrate the robustness and efficiency of
  the stabilized LCL procedure.
\end{abstract}



\pagestyle{myheadings}
\thispagestyle{plain}
\markboth{M. P. FRIEDLANDER AND M. A. SAUNDERS}%
         {A GLOBALLY CONVERGENT LCL METHOD}

\section{Introduction}                         \label{sec:introduction}

For optimization problems with nonlinear constraints, \emph{linearly
  constrained Lagrangian} (LCL) methods sequentially minimize a
Lagrangian function subject to linearized constraints.  As currently
defined, these methods converge rapidly near a solution but may not be
reliable from arbitrary starting points.  The well known example
\MINOS~\cite{MS82} has proven effective on many large and small
problems, especially within the \GAMS~\cite{BKM88}\ and
\AMPL~\cite{FGK93}\ environments, and is widely used in industry and
academia.  Its success motivates us to propose a globally convergent
variant of the LCL method.

Our stabilized LCL algorithm solves a sequence of linearly constrained
subproblems.  Each subproblem minimizes an augmented Lagrangian
function within a linear manifold that describes a current
approximation to the nonlinear constraints.  This manifold is
nominally a linearization of the constraint space but may be a relaxed
(i.e., larger) space at any stage, particularly during early
iterations.  Few conditions are imposed on the nature of the
subproblem solutions; consequently, the subproblems may be solved with
any of a variety of optimization routines for linearly constrained
problems, providing much flexibility.

The stabilized LCL method possesses two important properties: the
subproblems are always feasible, and they may be solved inexactly.
These features are present in \MINOS\ only as heuristics.  The method
may be regarded as a generalization of sequential augmented Lagrangian
methods (see, for example, \cite{GMW81,Ber82,Fle87}).  The theory we
develop provides a framework that unifies Robinson's LCL
method~\cite{Rob72} with the bound-constrained Lagrangian (BCL) method
used, for example, by \LANCE~\cite{CGT91a}.  In the context of our
theory, the proposed algorithm is actually a continuum of methods,
with LCL and BCL methods at opposite ends of a spectrum.  The
stabilized LCL algorithm exploits this connection between BCL and LCL
methods, preserving the fast local convergence properties of LCL
methods while inheriting the global convergence properties of BCL
methods.  This connection is explored in more detail by
Friedlander~\cite{Fri02}.

Our focus is on large-scale problems.  We implemented the stabilized
LCL method using the reduced-gradient part of \MINOS~\cite{MS78} and
the sequential quadratic programming code \SNOPT~\cite{GMS02} to solve
the linearly constrained subproblems.  These solvers are most
efficient on problems with few degrees of freedom.  Also, they use
only first derivatives, and consequently our implementation requires
only first derivatives.  We discuss how the stabilized LCL method
might be used with first- or second-derivative linearly constrained
solvers.

\subsection{The optimization problem}             \label{sec:intro-nlp}

The proposed method solves nonlinearly constrained optimization
problems of the form
\begin{equation*}
  \problem{\NP}{x\in\Real^n}
  {f(x)}
  {l \le
    \begin{pmatrix}
      x
\\    c(x)
\\    A x
    \end{pmatrix}
    \le u,}
\end{equation*}
where $f:\Real^n \mapsto \Real$ is a linear or nonlinear objective
function, $c:\Real^n \mapsto \Real^m$ is a vector of nonlinear
constraint functions, $A$ is a matrix, and $l$ and $u$ are vectors of
bounds.  We assume that $A$ and the derivatives of $c$ are sparse and
that the problem \NP\ is feasible.  We recognize that not all
optimization problems are feasible.  This possibility is addressed in
\Sec\ref{sec:sLCL-infeas}, where we explain how the proposed algorithm
reveals an infeasible optimization problem and discuss properties of
the points to which it converges.

One of the strengths of our method is that it does not explicitly
require second-order information.  However, the fast convergence rate
of the algorithm relies on sufficient smoothness of the nonlinear
functions, indicated by the existence of second derivatives. We make
that assumption:
\begin{assumption}
  \label{ass:continuity}
  The functions $f$ and $c$ are twice continuously differentiable on
  an open neighborhood containing the region
  \begin{equation*}
    l \le
    \begin{pmatrix}
      x \\ A x
    \end{pmatrix}
    \le u.
  \end{equation*}
\end{assumption}%
Note that second derivatives could be used if they were available,
thus accelerating the solutions of the subproblems and changing the
properties of the solutions obtained by the algorithm.  We discuss
this further in \Sec\ref{sec:sLCL-second}.

\subsection{The LCL approach}                    \label{sec:intro-lcl}

The acronym LCL is new.  Methods belonging to this class typically
have been referred to in the optimization literature as sequential
linearized constraint (SLC) methods (\cf\ \cite{GMW81,NW99}).  The
term SLC was chosen for compatibility with the terms sequential
quadratic programming (SQP) and sequential linear programming (SLP).
Those methods also sequentially linearize the constraints.  The term
\emph{linearly constrained Lagrangian}, however, emphasizes that the
Lagrangian itself, and not an approximation, is used in the
subproblems.  Moreover, there is a useful relationship (which we
exploit) between LCL and BCL methods, and this is hinted at by the
nomenclature.

The first LCL methods were proposed independently in 1972.
Robinson~\cite{Rob72} and Rosen and Kreuser~\cite{RK72} describe
similar algorithms based on minimizing a sequence of Lagrangian
functions subject to linearized constraints.  Robinson is able to
prove that, under suitable conditions, the sequence of subproblem
solutions converges quadratically to a solution of \NP.  A strength of
this method is that efficient large-scale methods exist for the
solution of the linearly constrained subproblems formed at each
iteration.  Any suitable example of these subproblem solvers may be
called as a black box.

\subsection{Other work on stabilizing LCL methods}
\label{sec:intro-other-lcl}

Other approaches to stabilizing LCL algorithms include two-phase
methods proposed by Rosen~\cite{Ros78} and Van Der Hoek~\cite{Vdh82}.
In these approaches, a Phase~1 problem is formed by moving the
nonlinear constraints into the objective by means of a quadratic
penalty function.  The solution of the Phase~1 problem is used to
initialize Robinson's method (Phase~2).  With a sufficiently large
penalty parameter, the Phase~1 solution will yield a starting point
that allows Robinson's method to converge quickly to a solution.
These two-phase methods choose the penalty parameter arbitrarily,
however, and do not deal methodically with infeasible linearizations.

In 1981, Best \etal~\cite{BBRR81} describe a variant of the two-phase
method whereby the Phase~1 penalty parameter is gradually increased by
repeated return to the Phase~1 problem if the Phase~2 iterations are
not converging.  This two-phase method differs further from Rosen's
and Van Der Hoek's methods in that the Phase~2 iterations involve only
those equality constraints identified as active by the Phase~1
problem.  The authors are able to retain local quadratic convergence
of the Phase~2 LCL iterations while proving global convergence to a
stationary point.  A drawback of their method is that it requires a
fourth-order penalty term to ensure continuous second derivatives of
the penalty objective.  This requirement may introduce significant
numerical difficulty for the solution of the Phase~1 problem (though
probably a quadratic-penalty term would suffice in practice).

Both two-phase methods share the disadvantage that the Phase~1 penalty
problems need to be optimized over a larger subspace than the
subsequent LCL phase.  We seek a method that retains the linearized
constraints as part of the subproblem, in order to keep the number of
degrees of freedom small; and, as in Robinson's 1972 method, we allow
the subproblem to determine the final set of active constraints.

\subsection{The generic problem}              \label{sec:intro-generic}

For the theoretical development of a stabilized LCL method, we
consider a simplified, generic formulation of \NP\ and take the
optimization problem to be
\begin{equation*}
  \problem{\GNP}{x\in\Real^n}
  {f(x)}
  {\begin{aligned}[t]
      c(x)  &=   0 \\
        x   &\ge 0,
    \end{aligned}}
\end{equation*}
where $c:\Real^n \mapsto \Real^m$.  Section~\ref{sec:imp} returns to
the formulation \NP\ in its discussion of the implementation of the
stabilized LCL method.

We define the \emph{augmented Lagrangian function} corresponding to
\GNP\ as
\begin{equation}
  \label{eq:50}
  \AugL(x,y,\rho) = f(x) - y\T c(x) + \half \rho \twonorm{c(x)}^2,
\end{equation}
where $x$, the $m$-vector $y$, and the scalar $\rho$ are independent
variables.  Let $g(x)$ denote the gradient of the objective function
$f(x)$, and let $J(x)$ denote the Jacobian matrix of the constraint
vector $c(x)$.  Denote by $H(x)$ and $H_i(x)$ the Hessian matrices of
$f(x)$ and $\compi{c(x)}$, respectively, where $\compi{\cdot}$ refers
to the $i$th component of a vector.  Define
\begin{equation}
  \label{eq:7}
  \fomult(x,y,\rho) = y - \rho c(x).
\end{equation}
The derivatives of $\AugL$ with respect to $x$ may be written as
follows:
\begin{align}
  \label{eq:96}
     \grad{x} \AugL(x,y,\rho)
     &= g(x) - J(x)\T \fomult(x,y,\rho)
\\\label{eq:16}
     \Hess{xx}\AugL(x,y,\rho)
     &= H(x) - \sum_{i=1}^m \compi{\fomult(x,y,\rho)} H_i(x)
     + \rho J(x)\T J(x).
\end{align}

We assume that problem \GNP\ is feasible and has at least one point
$(\xstar,\ystar,\zstar)$ that satisfies the first-order
Karush-Kuhn-Tucker (KKT) optimality conditions.
\begin{definition}[First-Order Optimality Conditions]  A triple
  $(\xstar,\ystar,\zstar)$ is a first-order KKT point for \GNP\ if for
  any $\rho \ge 0$ all of the following hold:
  \begin{subequations}  \label{eq:KKTGNP}
  \begin{align}
    c(\xstar)                           &=   0       \label{eq:KKTGNP1}
 \\ \grad{x} \AugL(\xstar,\ystar,\rho)  &=   \zstar  \label{eq:KKTGNP2}
 \\ \min(\xstar,\zstar)                 &=   0.      \label{eq:KKTGNP3}
  \end{align}
\end{subequations}
\end{definition}%
Note that~\eqref{eq:KKTGNP3} implies
\begin{subequations}
  \label{eq:9}
  \begin{align}
    \xstar &\ge 0        \label{eq:6}
\\  \zstar &\ge 0,       \label{eq:10}
  \end{align}
\end{subequations}
so that $\xstar$ and $\zstar$ must be primal and dual feasible,
respectively.

Let $\eta_*>0$ and $\omega_*>0$ be specified as primal and dual
convergence tolerances.  We regard the point $(x,y,z)$ to be an
acceptable solution of \GNP\ if it satisfies~\eqref{eq:KKTGNP} to
within these tolerances.  Specifically, we identify $(x,y,z)$ as an
approximate solution of \GNP\ if
\begin{subequations}  \label{eq:KKTGNPb}
  \begin{align}
    \norm{c(x)}                  &\le  \eta_*        \label{eq:KKTGNP1b}
 \\ \grad{x} \AugL(x,y,\rho)     &=    z             \label{eq:KKTGNP2b}
 \\ \infnorm{\min(x,z)}          &\le  \omega_*.     \label{eq:KKTGNP3b}
  \end{align}
\end{subequations}
Note that~\eqref{eq:KKTGNP3b} relaxes the nonnegativity
conditions~\eqref{eq:9} by the same tolerance $\omega_*$.  In
practice, we might choose to relax~\eqref{eq:6} to $x \ge -\delta_*
e$, for some $\delta_*>0$.  However, we ignore this detail for now.

For theoretical purposes, we assume that \emph{strict complementarity}
and the \emph{second-order sufficiency conditions} hold at each
$(\xstar,\ystar,\zstar)$.  We define these conditions as follows.
\begin{definition}[Strict Complementarity]
  \label{def:strict-comp}
  The point $(\xstar,\ystar,\zstar)$ satisfies strict complementarity
  if it satisfies~\eqref{eq:KKTGNP} and $\max(\xstar,\zstar) > 0$.
\end{definition}
\begin{definition}[Second-Order Sufficiency]
  \label{def:second-order}
  The point $(\xstar,\ystar,\zstar)$ satisfies the second-order
  sufficiency conditions for \GNP\ if it satisfies~\eqref{eq:KKTGNP}
  and strict complementarity and if for any $\rho \ge 0$,
  \begin{equation}
    \label{eq:3}
    p\T \Hess{xx}\AugL(\xstar,\ystar,\rho) p > 0
  \end{equation}
  for all $p\ne0$ satisfying
  \begin{equation}
    \label{eq:28}
    \text{$J(\xstar)p=0$ and $\compj{p}=0$ for all $j$
      such that $\compj \xstar=0$}
  \end{equation}
(and $\compj \zstar > 0$).
\end{definition}
\begin{assumption}
  \label{ass:second-order}
  The point $(\xstar,\ystar,\zstar)$ satisfies the second-order
  sufficiency conditions for \GNP.
\end{assumption}

\subsection{The canonical LCL method}         \label{sec:canonical-LCL}

The software package \MINOS\ solves the nonlinearly constrained
problem \GNP\ by minimizing a sequence of augmented Lagrangian
functions subject to linearized constraints.  Define the
\emph{constraint linearization} at the point $x\k$ as
\begin{equation*}
  \clink x = c(x\k) + J(x\k)(x - x\k).
\end{equation*}
Algorithm~\ref{algo:cLCL} outlines what we regard to be a canonical
LCL method.  It forms the basis for the \MINOS\ algorithm and is based
on solving the linearly constrained subproblems
\begin{equation*}
  \label{eq:LCk}
  \problem{\LCk}{x}{\AugL(x,y\k,\rho\k)}
  { \begin{aligned}[t]
      \clink x &=   0 \\
      x        &\ge 0,
    \end{aligned}
    }
\end{equation*}
which are parameterized by the latest estimates $x\k$ and $y\k$, and a
fixed penalty parameter $\rho\k\equiv\rhobar$ (which may be set to
zero).  The linear constraints $\clink x=0$ are the linearization of
$c$ at the point $x\k$.

Empirically, a positive penalty parameter $\rhobar$ has proven a
helpful addition to Robinson's method, but for other problems it has
been ineffective.  A theoretical understanding of when and how to
modify the penalty term has been lacking.

\begin{algo}{bt}

  \caption{Canonical LCL}

  \BlankLine

  \Inputs{$x_0, y_0, z_0$}
  \Outputs{$\xstar, \ystar, \zstar$}

  \BlankLine
  
  \Comment{Initialize parameters}{
    Set the penalty parameter
    $\rhobar \ge 0$.  Set positive convergence tolerances
    $\omega_*,\eta_* \ll 1$\;}

  $k \assign 0$\;
  \Converged $\assign$ \False\;
  
  \Repeat{\Converged}
  {
 
    \Comment{Solve the LC subproblem}{
      
      Solve \LCk\ to obtain a point $(\xstark,\dystark,\zstark)$.
      If there is more than one such point, choose
      $(\xstark,\dystark,\zstark)$ closest in norm to
      $(x\k,0,z\k)$\;
      $\ystark \assign y\k + \dystark$\;}
    
    \Comment{Update solution estimates}{%
        $x\kp1 \assign  \xstark$,
        $y\kp1 \assign  \ystark$,
        $z\kp1 \assign  \zstark$\;
      }

    \Comment{Test convergence}{%
    \lIf%
    {$(x\kp1,y\kp1,z\kp1)$ satisfies~\eqref{eq:KKTGNPb}}
    {\Converged $\assign$ \True}\;}
    \makebox[10ex][l]{%
    $\rho\k \assign \rhobar$;}
    [keep $\rho\k$ fixed]\\
    $k \assign k + 1$\;
  }
  $\xstar \assign x\k$, $\ystar \assign y\k$, $\zstar \assign z\k$\;
  \KwRet{$\xstar,\ystar,\zstar$}\;
  
  \label{algo:cLCL}
\end{algo}

\subsection{Notation}                        \label{sec:intro-notation}

The symbol $\xstar$ is used in two senses: as a limit point of the
sequence $\seq{x\k}$, and as the primal solution of \GNP.  We
distinguish between the two cases when the context is not clear.
Denote by $\ghat(x)$ the vector of components of $g(x)$ corresponding
to inactive bounds at $\xstar$, so that if $\Iscr =
\{i\in\iset{1}{n}\mid\compi{\xstar} > 0\}$\label{inactivebnds},
$\ghat(x) = \comp{\Iscr}{g(x)}$ (where $\comp{\Iscr}{\cdot}$ is a
shorthand notation for a subvector formed from the indices in
$\Iscr$).  Similarly, let $\Jhat(x)$ denote the corresponding columns
of the Jacobian matrix.

Unless otherwise specified, the function $\norm{x}$ represents the
Euclidean norm of the vector $x$.  When the arguments are vectors,
define the function $\min(\cdot,\cdot)$ component-wise.  The following
notation is used throughout:
\begin{center}
  \begin{tabular}{@{}ll@{}}
    $(x, y, z)$ & primal variables, dual variables, and reduced costs for \GNP,
\\  $(\xstar, \ystar, \zstar)$ & optimal variables for \GNP,
\\  $(x\k,y\k,z\k)$ & the $k$th estimate of $(\xstar,\ystar,\zstar)$,
\\  $(\xstark,\dystark,\zstark)$ & solution of the $k$th subproblem,
\\  $\ystark$ & $y\k + \dystark$; an updated multiplier estimate,
\\  $f\k$, $g\k$, $c\k$, $J\k$ & functions and gradients evaluated at $x\k$,
\\  $f_*$, $g_*$, $c_*$, $J_*$ & functions and gradients
                                         evaluated at $\xstar$.
  \end{tabular}
\end{center}
The augmented Lagrangian function is particularly important for our
analysis.  We often use the shorthand notation
\begin{equation}
  \label{eq:100}
  \AugL\k(x) \equiv \AugL(x,y\k,\rho\k)
             = f(x) - y\k\T c(x) + \half \rho\k \norm{c(x)}^2,
\end{equation}
when $y\k$ and $\rho\k$ are fixed.  The algorithms we discuss are
structured around \emph{major} and \emph{minor} iterations.  Each
major iteration solves a subproblem and generates an element of the
sequence $\seq{(x\k,y\k,z\k)}$.  Under certain (desirable)
circumstances, this sequence converges to a solution
$(\xstar,\ystar,\zstar)$.  For each major iteration $k$, there is a
corresponding set of minor iterations converging to
$(\xstark,\dystark,\zstark)$, the solution of the current subproblem.
In our development and analysis of a stabilized LCL method, we are
primarily concerned with the ``outer''-level algorithm.  Unless stated
otherwise, ``iterations'' refers to major iterations.

\section{An Elastic LC Subproblem}                  \label{sec:sLCL-LC}

The original LCL method introduced by Robinson~\cite{Rob72} sets
$\rhobar=0$ in Algorithm~\ref{algo:cLCL}.  A positive penalty
parameter could be used (as it is in \MINOS~\cite{MS82}) and may help
convergence from difficult starting points.

We recognize two particular causes of failure for the LCL method:
\begin{itemize}
\item The linearized constraints may be infeasible, so that the LCL
  iterations are not defined;
\item A near-singular Jacobian $J\k$ (we only assume nonsingularity of
  the Jacobian at limit points---\cf\ Assumption~\ref{ass:licq}) might
  lead to an arbitrarily large value of $\norm{\xstark-x\k}$
  regardless of the values of $y\k$ and $\rho\k$ in the subproblem
  objective.
\end{itemize}

To remedy both deficiencies we modify the linearized constraints used
by the LCL method, allowing some degree of flexibility in their
satisfaction.  We introduce a set of nonnegative elastic variables,
$v$ and $w$, into the constraints and introduce a penalty on these
variables into the subproblem objective.  Hence, we define the
subproblem as
\begin{equation*}
  \label{eq:ELCk}
  \problem{\ELCk}{x,v,w}
  {\AugL\k(x)
    + \sigma\k e\T    (v+w)}
  {\begin{aligned}[t]
      \clink x + v - w &=    0 \\ 
      x, v, w          &\geq 0,
    \end{aligned}}
\end{equation*}
where $e$ is a vector of ones.  \emph{This elastic subproblem is
  always feasible}.
Its solution yields a 5-tuple
$(\xstark,\dystark,\zstark,\vstark,\wstark)$ that satisfies the
first-order KKT conditions
\begin{subequations}
  \label{eq:58}
  \begin{align} 
   v,w                                  &\ge 0         \label{eq:58-1}
\\ \clink x + v - w                    &=   0         \label{eq:58-2}
\\ \grad{}\AugL\k(x) - J\k\T\Deltait y &=   z         \label{eq:58-3}
\\ \min(x,z)                           &=   0         \label{eq:58-4}
\\ \infnorm{\Deltait y}                &\le \sigma\k. \label{eq:58-5}
  \end{align}
\end{subequations}
Note that $\grad{}\AugL\k(x)$ involves $y\k$ and $\rho\k$.

The term $\sigma\k e\T (v+w)$ is the $\ell_1$-penalty function, and
together with the nonnegativity constraints $v,w\ge0$ it is equivalent
to a penalty on the one-norm of $(v - w)$.  We find later that the
bound~\eqref{eq:58-5} is crucial for the global convergence analysis
of the proposed method.

We note that the elastic LC subproblem can be equivalently stated as
\begin{equation*}
  \problem{\ELCkp}{x}
  { \AugL\k(x) + \sigma\k \onenorm{\clink x} }
  { x \ge 0, }
\end{equation*}
with solution $(\xstark,\zstark)$.  This immediately reveals the
stabilized LCL method's intimate connection with both the augmented
Lagrangian function and the BCL method.  Far from a solution, the
$\ell_1$-penalty term $\sigma\k \onenorm{\clink x}$ 
gives the method an opportunity to deviate from the constraint linearizations.
Near a solution, it keeps the iterates close to the linearizations.  
For values of $\sigma\k$ over a threshold value, the linearized
constraints are satisfied exactly, as required by the LCL method.

\subsection{The $\ell_1$-penalty function}       \label{sec:sLCL-ell1}

For any given subproblem of the stabilized LCL method, the penalty
term $\sigma\k e\T(v+w)$ may or may not equal zero, indicating that
the linearized constraints may not always be satisfied.  In contrast,
the \MINOS\ or the canonical LCL subproblems must always satisfy the
linearized constraints.  Thus, the set of active linearized
constraints of the stabilized LCL subproblem is always a subset
(though not necessarily strict) of the the canonical LCL subproblem.
Fletcher~\cite{Fle84}\ makes the same observation in connection with
his \SLQP\ method.  The global convergence properties of the
stabilized LCL method do not require independent constraint gradients
or bounded multipliers for each subproblem (these are required only at
limit points of the sequence generated by the algorithm).

\subsubsection*{Recovering the BCL subproblem}

Set $\sigma\k = 0$.  Then \ELCk\ and \ELCkp\ reduce to the
equivalent bound-constrained minimization problem
\begin{equation*}
  \label{eq:BCk}
  \problem{\BCk}{x}{\AugL\k(x)}{x \ge 0,}
\end{equation*}
where the bounds on the variables $v$ and $w$ have been eliminated
because they no longer appear in the objective.  The subproblem \BCk\ 
is used by the BCL method (see, for example, Hestenes~\cite{Hes69},
Powell~\cite{Pow69}, Bertsekas~\cite{Ber82}, and Conn
\etal~\cite{CGST96,CGT91b}).

\subsubsection*{Recovering the LCL subproblem}

The $\ell_1$-penalty function is \emph{exact}.  If the linearization
is feasible and $\sigma\k$ is larger than a certain threshold, $v$ and
$w$ will be zero and the minimizers of the elastic problem \ELCk\ will
coincide with the minimizers of the inelastic problem \LCk.  Exact
penalty functions have been studied by \cite{Ber82,Fle84,Lue84} among
others.  See the book by Conn \etal~\cite{CGT00} for a more recent
discussion.

We are particularly interested in this feature when the iterates
generated by the stabilized LCL algorithm are approaching a solution
$(\xstar,\ystar,\zstar)$.  Recovering the canonical LCL subproblem as
the iterates approach a solution ensures that the stabilized LCL
method inherits LCL's fast local convergence properties.

To prove that the condition
\begin{equation}
  \label{eq:4}
  \infnorm{\dystark} < \sigma\k
\end{equation}
is sufficient to force the elastic variables to zero, we require two
conditions: (i) the inelastic subproblem \LCk\ must satisfy the
second-order sufficiency conditions at a solution $\xstark$; and (ii)
$\xstark$ must be a regular point.  Assumptions~\ref{ass:second-order}
and~\ref{ass:licq} guarantee that both these conditions are met.  For
$\xstark$ close to $\xstar$, Assumption~\ref{ass:second-order}
guarantees that \LCk\ satisfies the second-order conditions.
Assumption~\ref{ass:licq} guarantees the regularity of $\xstark$ when
it is near $\xstar$.  Lemma~\ref{le:L1threshold} establishes the
threshold value of $\sigma\k$.

\begin{lemma}
  \label{le:L1threshold}
  Suppose that $(\xstark,\dystark,\zstark)$ satisfies the second-order
  sufficiency conditions for \LCk. Then if~\eqref{eq:4} holds,
  $(\xstark,\dystark,\zstark)$ also solves \ELCk.
\end{lemma}

\begin{proof}
  See Luenberger~\cite[p.~389]{Lue84}.
\end{proof}

\subsection{Early termination of the subproblems}\label{sec:sLCL-earlyterm}

Poor values of $x\k$, $y\k$, or $\rho\k$ may imply subproblems whose
accurate solutions are far from a solution of \GNP.  We therefore
terminate subproblems early by relaxing~\eqref{eq:58-4}
and~\eqref{eq:58-5} by an amount $\omega\k$.  However, we enforce the
nonnegativity condition on $x$ (implied by~\eqref{eq:58-4}):
\begin{subequations}
  \label{eq:59}
\begin{align}
   x,v,w                                 &\ge 0           \label{eq:59-a}
\\ \clink x + v - w                     &=   0           \label{eq:59-3}
\\ \grad{}\AugL\k(x) - J\k\T \Deltait y &=   z           \label{eq:59-c}
\\ \infnorm{\min(x,z)}                  &\le \omega\k    \label{eq:59-d}
\\ \infnorm{\Deltait y}                 &\le \sigma\k + \omega\k.   \label{eq:59-6}
\end{align}
\end{subequations}
Each subproblem is required to return a solution satisfying the linear
and nonnegativity constraints and, as discussed in connection
with~\eqref{eq:6}, in practice \eqref{eq:59-a} and/or~\eqref{eq:59-3}
would be relaxed by a fixed tolerance $\delta$.

\section{The Stabilized LCL Algorithm}            \label{sec:sLCL-algo}

Algorithm~\ref{algo:sLCL} outlines the stabilized LCL method.  Its
structure closely parallels the BCL algorithm described
in~\cite{CGT91b}.  Based on the current primal infeasibility, each
iteration of the algorithm is regarded as either ``successful'' or
``unsuccessful.''  In the ``successful'' case, the solution estimates
are updated by using information from the current subproblem solution.
If the iteration is ``unsuccessful,'' the subproblem solutions are
discarded, the current solution estimates are held fixed, and the
penalty parameter $\rho\k$ is increased in an effort to reduce the
primal infeasibility in the next iteration.  In order for the
linearized constraints not to continue interfering with the penalty
parameter's ability to reduce the primal infeasibility, the algorithm
relaxes the linearizations by reducing the elastic penalty parameter
$\sigma\k$.

The two salient features of this algorithm are that it is globally
convergent and that it is asymptotically equivalent to the canonical
LCL method.  In \Sec\ref{sec:sLCL-global} we demonstrate the global
convergence properties of the algorithm by proving results analogous
to Lemma~4.3 and Theorem~4.4 in~\cite{CGT91b}.  In
\Sec\ref{sec:sLCL-local} we demonstrate that the algorithm eventually
reduces to the canonical LCL method and hence inherits that method's
asymptotic convergence properties.

\begin{algo}{tbp} \small
  
  \caption[The stabilized LCL algorithm]{Stabilized LCL.}

  \BlankLine

  \Inputs{$x_0, y_0, z_0$}
  \Outputs{$\xstar, \ystar, \zstar$}

  \BlankLine
  
  \Comment{Initialize parameters}{
    Set $\sigmaB > \sigmab > 0$.  Set constants $\tau_\rho,\tau_\sigma
    > 1$.  Set the initial penalty parameters $\rho_0 > 1$ and
    $\sigma_0 \gg 1$.  Set positive convergence tolerances
    $\omega_*,\eta_* \ll 1$ and initial tolerances $\omega_0 >
    \omega_*$ and $\eta_0 > \eta_*$. Set constants $\alpha, \beta>0$
    with $\alpha<1$\;}

  $k \assign 0$\;
  \Converged $\assign$ \False;

  \Repeat{\Converged}
  {
    \lnl{algo:sLCLomega}
    Choose $\omega\k \ge \omega_*$ such that
    $\lim_{k\to\infty}\omega\k = \omega_*$\;
    
    \lnl{algo:sLCLsolveELC}
    \Comment{Solve the LC subproblem}{
      Solve \ELCk\ to obtain a point $(\xstark,\dystark,\zstark)$
      satisfying~\eqref{eq:59}.  If there is more than one such point,
      compute the one closest in norm to $(x\k,0,z\k)$\;
      $\ystark \assign y\k + \dystark$\;}
    \lnl{algo:sLCLbranchtest}
    \eIf{$\norm{c(\xstark)} \leq \max(\eta_*, \eta\k)$}
    {
      \lnl{algo:sLCLupdatesolutions}
      \Comment{Update solution estimates}{
        $x\kp1 \assign \xstark$\;
        \lnl{algo:sLCLmult}
        \makebox[41ex][l]{%
          $y\kp1 \assign \ystark - \rho\k c(\xstark)$ $(\equiv
          \fomult(\xstark,\ystark,\rho\k))$;}
        [or $y\kp1 \assign \ystark$]\\
        \makebox[41ex][l]{%
        $z\kp1 \assign \zstark$;}
        [or $z\kp1 \assign g\kp1 - J\kp1\T y\kp1$]
      }

      \lnl{algo:sLCLupdatesigma}
      \Comment{Update penalty parameter and elastic weight}{
        \makebox[41ex][l]{%
        $\rho\kp1   \assign  \rho\k$;}
        [keep $\rho\k$]\\
        \makebox[41ex][l]{%
        $\sigma\kp1 \assign
        \max\{\sigmab,
        \min(\infnorm{\dystark}, \sigmaB)\}$;}
        [reset $\sigma\k$]
      }
      
      \lnl{algo:sLCLtestconv}
      \Comment{Test convergence}{
        \lIf{$(x\kp1,y\kp1,z\kp1)$ satisfies~\eqref{eq:KKTGNPb}}
        {\Converged $\assign$ \True\;}}
      \lnl{algo:sLCLetab}
      \makebox[45ex][l]{%
      $\eta\kp1 \assign \eta\k / \rho\kp1^{\beta}$\;}
      [decrease $\eta\k$]
    } 
    {
      \lnl{algo:sLCLkeep}
      \Comment{Keep solution estimates}{
      $x\kp1 \assign x\k$;
      $y\kp1 \assign y\k$;
      $z\kp1 \assign z\k$\;}
      
      \lnl{algo:sLCLincrho}
      \Comment{Update penalty parameter and elastic weight}{
        \makebox[41ex][l]{%
        $\rho\kp1   \assign \tau_\rho \rho\k$;}
        [increase $\rho\k$]\\
        \makebox[41ex][l]{%
        $\sigma\kp1 \assign \sigma\k / \tau_\sigma$;}
       [decrease $\sigma\k$]
     }
      \lnl{algo:sLCLetaa}
      \makebox[45ex][l]{%
        $\eta\kp1 \assign \eta_0 / \rho\kp1^{\alpha}$\;}
      [may increase or decrease $\eta\k$]
    } 
    $k \assign k + 1$\;
  } 
  $\xstar \assign x\k$;
  $\ystar \assign y\k$;
  $\zstar \assign z\k$\;
  \KwRet{$\xstar, \ystar, \zstar$}\;

  \label{algo:sLCL}

\end{algo}

\clearpage
\subsection{Global convergence properties}      \label{sec:sLCL-global}

We make the following assumptions.

\begin{assumption}
  \label{ass:compactness}
  The sequence of iterates $\seq\xstark$ lies in the closed and
  bounded set $\Bscr \subset \Real^n$.
\end{assumption}
\begin{assumption}
  \label{ass:licq}
  The matrix $\Jhat(\xstar)$ has full row rank at every limit point
  $\xstar$ of the sequence $\seq\xstark$.
\end{assumption}

The first assumption guarantees that any sequence of iterates
generated by the algorithm always has some convergent subsequence.
The second assumption is commonly known as the \emph{linear
  independence constraint qualification} (LICQ) (see, for example,
Mangasarian~\cite{Man69}, or for a more recent reference, Nocedal and
Wright~\cite{NW99}).

Let $\xstar$ be any limit point of the sequence $\seq\xstark$.  At all
points $x$ for which $\Jhat(x)$ has full row rank we define the
\emph{least-squares multiplier estimate}, $\lsmult(x)$, as the
solution of the linear least-squares problem
\begin{equation}
  \label{eq:2}
  \minimize{y} \ \norm{\ghat(x) - \Jhat(x)\T y}^2.
\end{equation}
Note that the definitions of $\ghat$, $\Jhat$, and hence $\lsmult$
require a priori knowledge of the bounds active at $\xstar$.  We
emphasize that $\lsmult$ is used only as an analytical device and its
computation is never required.  Assumption~\ref{ass:licq} guarantees
the uniqueness of $\lsmult$ at every limit point of the sequence $\seq
\xstark$.

\subsubsection{Convergence of LC subproblem solutions}

In this section we prove that the sequence of LC subproblem solutions
generated by Algorithm~\ref{algo:sLCL} converges to a KKT point of
\GNP.

We need the following lemma to bound the errors in the least-squares
multiplier estimates relative to the error in $x\k$.  The lemma simply
demonstrates that $\lsmult(x)$ is Lipschitz continuous in a
neighborhood of $\xstar$.

\begin{lemma} \label{le:lambda}
  Let $\seq{x\k}$, $k\in\Kscr$ be a subsequence converging to $\xstar$
  and suppose that Assumptions~\aref{ass:continuity}
  and~\aref{ass:licq} hold.  Then there exists a positive constant
  $\alpha$ such that $\norm{\lsmult(x\k) - \lsmult(\xstar)} \leq
  \alpha \norm{x\k - \xstar}$ for all $k\in\Kscr$ sufficiently large.
\end{lemma}

\begin{proof}
  See Lemmas~2.1 and~4.4 of~\cite{CGST96}.
\end{proof}

To prove the global convergence properties of
Algorithm~\ref{algo:sLCL}, we first describe the properties of any
limit point that the algorithm generates.  We are not claiming (yet!)
that the algorithm is globally convergent, only that if it \emph{does}
converge, then the set of limit points generated must satisfy some
desirable properties.  The following lemma is adapted from Lemma 4.4
of~\cite{CGST96}.

\begin{lemma} \label{le:penalty}
  Let $\seq{\omega\k}$ and $\seq{\rho\k}$ be sequences of positive
  scalars, where $\omega\k\to0$.  Let $\seq{x\k}$ be any sequence of
  $n$-vectors and $\seq{y\k}$ be any sequence of $m$-vectors.  Let
  $\seq{(\xstark,\dystark,\zstark)}$ be a sequence of vectors
  satisfying~\eqref{eq:59-a}, \eqref{eq:59-c}, and \eqref{eq:59-d}.
  Let $\xstar$ be any limit point of the sequence $\seq\xstark$, and
  let $\Kscr$ be the infinite set of indices associated with that
  convergent subsequence.  Suppose that
  Assumptions~\aref{ass:continuity}, \aref{ass:compactness},
  and~\aref{ass:licq} hold.  Set $\ystark=y\k+\dystark$, $\fomult\k =
  \fomult(\xstark,\ystark,\rho\k)$, and $\ystar = \lsmult(\xstar)$.
  The following properties then hold:
  \begin{enumerate}
  \item \label{le:penalty:1} There are positive constants $\alpha_1$,
    $\alpha_2$, and $M$ such that
    \begin{align}
      \label{eq:22}
      \norm{\fomult\k - \ystar}
      &\leq \beta_1 \equiv
      \alpha_1 \omega\k + M \norm{\xstark - x\k} \,
      \norm{\ystark - y\k} +
      \alpha_2 \norm{\xstark - \xstar},
      \\
      \label{eq:cconv1}
      \rho\k\norm{c(\xstark)}
      &\leq
      \beta_2 \equiv \beta_1 +
      \norm{\ystark - y\k} + \norm{y\k - \ystar},
    \end{align}
    for all $k\in\Kscr$ sufficiently large.
  \item \label{le:penalty:2} As $k\in\Kscr$ gets large, if
    $\norm{\ystark - y\k} \to 0$, or if $\norm{\ystark - y\k}$ is
    bounded and $\norm{\xstark - x\k} \to 0$, then
    \begin{equation*}
      \fomult\k \to \ystar
      \text{and}
      \zstark   \to \zstar \defd \grad{x}\AugL(\xstar,\ystar,0).
    \end{equation*}
  \item \label{le:penalty:3} If, in addition, $\cstar = 0$, then
    $(\xstar,\ystar,\zstar)$ is a first-order KK point for \GNP.
  \end{enumerate}
\end{lemma}

\begin{proof}
  From the definition of $\lsmult(\xstark)$, the least-squares
  multiplier estimates,
  \begin{equation}
    \label{eq:63}
    \begin{aligned}
      \norm{\lsmult(\xstark) - \fomult\k}
      &=  \norm{(\Jhat(\xstark) \Jhat(\xstark)^T)\inv
        \Jhat(\xstark) \ghat(\xstark) - \fomult\k} \\
      &=  \norm{
        ( \Jhat(\xstark) \Jhat(\xstark)^T )\inv \Jhat(\xstark)
        ( \ghat(\xstark) - \Jhat(\xstark)\T \fomult\k )
        } \\
      &\leq \norm{(\Jhat(\xstark) \Jhat(\xstark)^T)\inv \Jhat(\xstark)} 
        \cdot \norm{\ghat(\xstark) - \Jhat(\xstark)\T \fomult\k}.
    \end{aligned}
  \end{equation}
  By assumption, $\Jhat(\xstar)$ has full row rank.  Continuity of $J$
  then implies that
  \begin{equation*}
    (\Jhat(\xstark) \Jhat(\xstark)^T)\inv \Jhat(\xstark)
  \end{equation*}
  exists for all $k\in\Kscr$ large enough.  Then
  there exists a positive scalar $\alpha_1$ such that
  \begin{equation}
    \label{eq:64}
    \norm{(\Jhat(\xstark) \Jhat(\xstark)^T)\inv
      \Jhat(\xstark)} \leq \frac{\alpha_1}{\sqrt n},
  \end{equation}
  where $n$ is the dimension of the vector $x$.
  Substituting~\eqref{eq:64} into~\eqref{eq:63},
  \begin{equation}
    \label{eq:65}
    \norm{\lsmult(\xstark) - \fomult\k}
    \leq \frac{\alpha_1}{\sqrt n}
    \norm{\ghat(\xstark) - \Jhat(\xstark)\T \fomult\k}.
  \end{equation}
  We now show that $\norm{\ghat(\xstark) - \Jhat(\xstark)\T\fomult\k}$
  is bounded.  By hypothesis, $(\xstark,\ystark,\zstark)$
  satisfies~\eqref{eq:59-c}.  Together with~\eqref{eq:96},
  \begin{equation}
    \label{eq:40}
    \begin{aligned} 
      \zstark
        &= \grad{}\AugL\k(\xstark) - J\k\T\dystark
      \\&= g(\xstark) - J(\xstark)\T (y\k - \rho\k c(\xstark))
           - J\k\T \dystark
      \\&= g(\xstark)
           - J(\xstark)\T  \bigl(y\k+\dystark-\rho\k c(\xstark) \bigr) 
           + \big(J(\xstark) - J\k\big)\T \dystark
      \\&= g(\xstark)
           - J(\xstark)\T  \fomult\k
           + \big(J(\xstark) - J\k\big)\T (\ystark - y\k),
  \end{aligned}
  \end{equation}
  where $\ystark\defd y\k+\dystark$ and $\fomult\k \defd
  \fomult(\xstark,\ystark,\rho\k)=\ystark-\rho\k c(\xstark)$.  For
  $k\in\Kscr$ large enough, $\xstark$ is sufficiently close to
  $\xstar$ so that
  \begin{equation}
    \label{eq:1}
    \norm{\comp{\Iscr}{\zstark}} \le \norm{\min(\xstark,\zstark)},
  \end{equation}
  where $\Iscr$ is the index set of inactive bounds at $\xstark$, as
  defined in \Sec\ref{inactivebnds}.  Because $\xstark$ and $\zstark$
  both satisfy~\eqref{eq:59-d}, \eqref{eq:1} implies that
  \begin{equation}
    \label{eq:18}
    \norm{\comp{\Iscr}{\zstark}} \le \sqrt n \ \omega\k.
  \end{equation}
  Combining~\eqref{eq:40} and~\eqref{eq:18},
  \begin{equation}
    \label{eq:37}
    \norm{\ghat(\xstark) - \Jhat(\xstark)\T \fomult\k
      + \big(\Jhat(\xstark) - \Jhat\k\big)\T (\ystark - y\k)}
    \le
    \sqrt n \ \omega\k.
  \end{equation}
  But, from the triangle and Cauchy-Schwartz inequalities, we have
  \begin{multline}
    \label{eq:66}
    \norm{\ghat(\xstark) - \Jhat(\xstark)\T \fomult\k}
    \leq
    \norm{\ghat(\xstark) - \Jhat(\xstark)\T \fomult\k +
      \big(\Jhat(\xstark) - \Jhat\k\big)\T (\ystark - y\k)} \\
    + \norm{\Jhat(\xstark) - \Jhat\k} \norm{\ystark - y\k}.
  \end{multline}
  Also, the continuity of $J$ implies that there exists a positive
  constant $M$ such that $\norm{\Jhat(\xstark) - \Jhat\k} \leq M
  \frac{\sqrt n}{\alpha_1}\norm{\xstark - x\k}$.  Together,
    \eqref{eq:66} and~\eqref{eq:37} imply that
  \begin{equation}
    \label{eq:67}
    \norm{\ghat(\xstark) - \Jhat(\xstark)\T \fomult\k} 
    \leq
    \sqrt n \ \omega\k + M\frac{\sqrt n}{\alpha_1}\norm{\xstark - x\k}
    \norm{\ystark - y\k},
  \end{equation}
  and so we have derived a bound on $\norm{\ghat(\xstark) -
    \Jhat(\xstark)\T\fomult\k}$, as required.
  
  We now derive~\eqref{eq:22}.  From the triangle inequality,
  \begin{equation}
    \label{eq:98}
    \norm{\fomult\k - \ystar}
    \leq \norm{\lsmult(\xstark) - \fomult\k}
    +    \norm{\lsmult(\xstark) - \ystar}.
  \end{equation}
  Using inequality~\eqref{eq:67} in~\eqref{eq:65}, we deduce that
  \begin{equation}
    \label{eq:68}
    \norm{\lsmult(\xstark) - \fomult\k}
    \leq
    \alpha_1 \omega\k + M\norm{\xstark - x\k}
    \norm{\ystark - y\k},
  \end{equation}
  and Lemma~\ref{le:lambda} implies that there exists a constant
  $\alpha_2$ such that
  \begin{equation}
    \label{eq:45}
    \norm{\lsmult(\xstark)-\ystar}\le\alpha_2\norm{\xstark-\xstar},
  \end{equation}
  for all $k\in\Kscr$ large enough (recall that
  $\ystar\equiv\lsmult(\xstar)$).  Substituting~\eqref{eq:68}
  and~\eqref{eq:45} into~\eqref{eq:98}, we obtain $\norm{\fomult\k -
    \ystar} \leq \beta_1$ as stated in~\eqref{eq:22}.
  
  We now prove~\eqref{eq:cconv1}.  From the definition of $\fomult\k$,
  rearranging terms yields
  \begin{equation}
    \label{eq:70}
    \rho\k c(\xstark) = \ystark - \fomult\k.
  \end{equation}
  Taking norms of both sides of~\eqref{eq:70} and using~\eqref{eq:22}
  yields
  \begin{equation*}
    \begin{aligned}
      \rho\k\norm{c(\xstark)}
      &=   \norm{\ystark - \fomult\k} \\
      &=   \norm{y\k - \ystar + \ystar - \fomult\k
         + \ystark - y\k}\\
      &\le \norm{\fomult\k - \ystar}
         + \norm{y\k - \ystar}
         + \norm{\ystark - y\k}\\
      &\le \beta_1
         + \norm{y\k - \ystar} + \norm{\ystark - y\k}
\\    &\equiv \beta_2,
    \end{aligned}
  \end{equation*}
  and so Part~1 of Lemma~\ref{le:penalty} is proved.
  
  Now suppose that $\norm{\ystark - y\k} \to 0$ as $k\in\Kscr$ goes to
  infinity.  Because $\seq\xstark$ and $\seq{x\k}$ are in the compact
  set $\Bscr$, $\norm{\xstark - x\k}$ is bounded.  We conclude
  from~\eqref{eq:22} that $\fomult\k \to \ystar$ as $k\in\Kscr$ goes
  to infinity.  We also conclude from the continuity of $J$ that
  $\norm{J(\xstark) - J\k}$ is bounded, so that
  \begin{equation}
    \label{eq:20}
    \lim_{k\in\Kscr} \norm{(J(\xstark) - J\k)\T(\ystark - y\k)} = 0.
  \end{equation}
  On the other hand, suppose that $\norm{\ystark - y\k}$ is uniformly
  bounded and that $\lim_{k\in\Kscr}\norm{\xstark - x\k}=0$.  We then
  conclude from~\eqref{eq:22} that $\fomult\k \to \ystar$ as
  $k\in\Kscr$ goes to infinity and~\eqref{eq:20} holds.  Because
  $\lim_{k\in\Kscr}(\xstark,\fomult\k)=(\xstar,\ystar)$,
  \begin{equation*}
    \label{eq:21}
    g(\xstark) - J(\xstark)\T\fomult\k \to
    \gstar     -\Jstar\T     \ystar,
  \end{equation*}
  and so~\eqref{eq:40} and~\eqref{eq:20} together imply that
  \begin{equation}
    \label{eq:39}
    \zstark \to \zstar \equiv \grad{x}\AugL(\xstar,\ystar,0)
  \end{equation}
  as $k\in\Kscr$ goes to infinity.  Thus we have proved Part~2 of
  Lemma~\ref{le:penalty}.

  Now suppose that
  \begin{equation}
    \label{eq:43}
      \cstar = 0.
  \end{equation}
  Each $\xstark$ and $\zstark$ satisfies~\eqref{eq:59-d}.  Then
  $\lim_{k\in\Kscr}(\xstark,\zstark) = (\xstar,\zstar)$, and $\omega\k
  \to 0$ implies
  \begin{equation}
    \label{eq:42}
    \min(\xstar,\zstar) = 0.
  \end{equation}
  Therefore,~\eqref{eq:39}--\eqref{eq:42} imply that
  $(\xstar,\ystar,\zstar)$ satisfies~\eqref{eq:KKTGNP} and so it is a
  first-order KKT point for \GNP.  Part~3 is thus proved, and the
  proof is complete.
\end{proof}

The conclusions of Lemma~\ref{le:penalty} pertain to any sequence
$\seq{(\xstark,\dystark,\zstark)}$ satisfying the approximate
first-order conditions~\eqref{eq:59}.  Algorithm~\ref{algo:sLCL}
generates such a sequence and also generates auxiliary sequences of
scalars $\seq{\omega\k}$, $\seq{\rho\k}$, and $\seq{\sigma\k}$ in such
a way as to guarantee that the hypotheses of Lemma~\ref{le:penalty}
hold.  We demonstrate in Theorem~\ref{th:globalconv} that the
condition of Part~3 of Lemma~\ref{le:penalty} holds.  Therefore, every
limit point of the sequence $\seq{(\xstark,\fomult\k,\zstark)}$ is a
first-order KKT point for \GNP.

\subsubsection{Convergence of $\norm{y\k}/\rho\k$}

Before laying out the global convergence properties of the stabilized
LCL method, we need to show that if $\rho\k \to \infty$
then the quotient $\norm{y\k}/\rho\k$
converges to 0.  This property is required (and used by Conn
\etal~\cite{CGST96,CGT91b}) in lieu of assuming that $\norm{y\k}$
remains bounded.

\begin{lemma}
  \label{le:lambdarhosLCL}
  Suppose that $\rho\k \to \infty$ as $k$ increases when
  Algorithm~\aref{algo:sLCL} is executed.  Then $\norm{y\k} / \rho\k
  \to 0$.
\end{lemma}

\begin{proof}
  The multiplier update of Algorithm~\ref{algo:sLCL}
  (Step~\ref{algo:sLCLmult}) is
  $y\kp1 \assign \ystark - \rho\k c(\xstark)$.
  If this is replaced by
  $y\kp1 \assign y\k - \rho\k c(\xstark)$,
  Lemma~4.1 of Conn \etal~\cite{CGT91b} applies. The construction of
  the forcing sequence $\eta\k$ (Steps~\ref{algo:sLCLetab}
  and~\ref{algo:sLCLetaa}) is therefore sufficient to guarantee that
  $\norm{y\k} / \rho\k \to 0$.  Note that the norm of the difference
  between the two updates
  is given by $\norm{\ystark -
    y\k}\equiv\norm{\dystark}$.  This difference is bounded because
  $\dystark$ satisfies~\eqref{eq:59-6}.  Therefore, $\norm{y\k} /
  \rho\k \to 0$ in Algorithm~\ref{algo:sLCL} as $\rho\k \to \infty$.
\end{proof}


\subsubsection{Main convergence result}

With Lemmas \ref{le:penalty}--\ref{le:lambdarhosLCL} in hand,
we are now able to prove global convergence of the stabilized LCL method.

\begin{theorem} \label{th:globalconv}
  Let $\seq{(\xstark,\ystark,\zstark)}$ be the sequence of vectors
  generated by Algorithm~\aref{algo:sLCL} with tolerances $\omega_*=0$
  and $\eta_*=0$.  Let $\xstar$ be any limit point of the sequence
  $\seq\xstark$ and let $\Kscr$ be the infinite set of indices
  associated with that convergent subsequence.  Then, under the
  assumptions of Lemma~\aref{le:penalty}, Parts~\aref{le:penalty:1},
  \aref{le:penalty:2}, and~\aref{le:penalty:3} of that lemma hold.
\end{theorem}
\begin{proof}
  Algorithm~\ref{algo:sLCL} generates positive scalars $\rho\k$ and,
  by Steps~\ref{algo:sLCLomega}, \ref{algo:sLCLetab},
  and~\ref{algo:sLCLetaa}, generates positive scalars $\omega\k \to 0$
  and $\eta\k \to 0$.  Step~\ref{algo:sLCLsolveELC} of the algorithm
  generates a sequence $\seq{(\xstark,\ystark,\zstark)}$, where
  $\ystark\equiv y\k+\dystark$.  Each $(\xstark,\dystark,\zstark)$
  satisfies~\eqref{eq:59}. Therefore, the hypotheses of
  Lemma~\ref{le:penalty} hold, and Part~1 of the lemma follows
  immediately.
  
  Note that each $\xstark$ satisfies~\eqref{eq:59-a}, and so
  $\xstark\ge0$ for all $k$.  Thus, $\xstar\ge0$.  Moreover, because
  $\tau_\sigma > 1$ and $\sigmaB$ is finite,
  Steps~\ref{algo:sLCLupdatesigma} and~\ref{algo:sLCLincrho} of
  Algorithm~\ref{algo:sLCL} ensure that $\sigma\k$ is uniformly
  bounded.  We then need to consider the four possible cases. For all
  $k\in\Kscr$,
  \begin{enumerate}
  \item 
    \label{ca:rhoBsig0}
    $\rho\k$ is uniformly bounded, and $\sigma\k \to 0$ as $k$ gets
    large;
  \item
    \label{ca:rhoBsigB}
    $\rho\k$ is uniformly bounded, and $\sigma\k$ is uniformly
    bounded away from zero;
  \item
    \label{ca:rhoIsig0}
    $\rho\k \to \infty$ and $\sigma\k \to 0$ as $k$ gets large;
  \item
    \label{ca:rhoIsigB}
    $\rho\k \to \infty$ and $\sigma\k$ is uniformly bounded away from
    zero.
  \end{enumerate}
  For the remainder of this proof, we consider only $k\in\Kscr$.

  We dismiss Case~\ref{ca:rhoBsig0} because it cannot be generated by
  the algorithm. (As $k$ gets large, $\sigma\k \to 0$ only if
  Step~\ref{algo:sLCLincrho} is executed infinitely many times,
  contradicting the finiteness of $\rho\k$.)
  
  Case~\ref{ca:rhoBsigB} implies that Step~\ref{algo:sLCLupdatesigma}
  of Algorithm~\ref{algo:sLCL} is executed for all $k$ large enough.
  Thus, $x\kp1=\xstark$ for all large $k$, and hence $\xstark \to
  \xstar$ implies $x\k \to \xstar$.  Therefore, $\norm{\xstark - x\k}
  \to 0$.  Because each $\dystark$ satisfies~\eqref{eq:59-6},
  $\ystark$ satisfies
  \begin{equation}
    \label{eq:23}
    \infnorm{\ystark - y\k} \le \omega\k + \sigma\k.
  \end{equation}
  Because $\sigma\k$ and $\omega\k$ are uniformly bounded, Part~2 of
  Lemma~\ref{le:penalty} holds.  In addition,
  $\norm{c(\xstark)}\le\eta\k$ for all $k$ large enough, and so
  $\eta\k \to 0$ implies that $c(\xstark) \to 0$.  By continuity of
  $c$, $\cstar = 0$.  Thus, Part~3 of Lemma~\ref{le:penalty} holds.
  
  Now consider Case~\ref{ca:rhoIsig0}.  Because $\sigma\k \to 0$ and
  $\omega\k \to 0$, \eqref{eq:23} implies that $\norm{\ystark - y\k}
  \to 0$ as $k$ increases.  Then Part~2 of the lemma holds.  To show
  that $c(\xstark) \to 0$, divide both sides of~\eqref{eq:cconv1} by
  $\rho\k$ to obtain
  \begin{equation*}
    \norm{c(\xstark)}
    \leq
    \underbrace{\frac{\alpha_1\omega\k}{\rho\k}}_{(a)}
    +
    \underbrace{
      \frac{1}{\rho\k}
      \norm{\ystark - y\k} \big(M \norm{\xstark - x\k} + 1\big)}_{(b)}
    +
    \underbrace{\frac{\alpha_2}{\rho\k} \norm{\xstark - \xstar}}_{(c)}
    +
    \underbrace{\frac{1}{\rho\k}\norm{y\k - \ystar}}_{(d)}.
  \end{equation*}
  Term $(a)$ clearly goes to zero as $\rho\k$ increases.  Because
  $\ystark$ and $y\k$ satisfy~\eqref{eq:23}, and because $\xstark$ and
  $x\k$ belong to the compact set $\Bscr$, $(b)$ and $(c)$ go to zero as
  $\rho\k$ increases.  By Lemma~\ref{le:lambdarhosLCL},
  $\norm{y\k}/\rho\k\to 0$, and so $(d)$ goes to 0.  We conclude that
  $\norm{c(\xstark)}\to 0$ as $k$ increases, as required.
  
  In Case~\ref{ca:rhoIsigB}, both Steps~\ref{algo:sLCLupdatesigma}
  and~\ref{algo:sLCLincrho} are executed infinitely often.  But
  because $\sigma\k$ is uniformly bounded, so is $\norm{\ystark -
    y\k}$.  As in Case~\ref{ca:rhoBsigB}, $\norm{\xstark-x\k} \to 0$
  as $k$ get large, and so Part~2 of Lemma~\ref{le:penalty} holds.
  The rest of the analysis for this case is the same as for
  Case~\ref{ca:rhoIsig0}.
\end{proof}

\subsubsection{Finite termination}

Note that the convergence test takes place only if
Step~\ref{algo:sLCLbranchtest} of Algorithm~\ref{algo:sLCL} tests
true; \ie, if $\norm{c(\xstark)} \le \eta\k$ (because $\eta_*=0$).  To
guarantee that the algorithm will eventually terminate as the iterates
$x\k$, $y\k$, and $z\k$ converge, we need to guarantee that
Steps~\ref{algo:sLCLupdatesolutions} and~\ref{algo:sLCLtestconv}
execute infinitely often.  The forcing sequence $\eta\k$ is intimately
tied to this occurrence.  For example, if $\eta\k \equiv 0$, then we
would not normally expect Step~\ref{algo:sLCLbranchtest} to evaluate
true (except in rare occasions when $c(\xstark) = 0$).  The forcing
sequence defined by Steps~\ref{algo:sLCLetab} and~\ref{algo:sLCLetaa}
of Algorithm~\ref{algo:sLCL} is suggested by Conn
\etal~\cite{CGST96,CGT91b}.  The following corollaries show that this
forcing sequence has the desired property and summarize the global
convergence properties of Algorithm~\ref{algo:sLCL}.  Unlike for the
previous results in this section, we now need to strengthen our
assumptions and require that only a single limit point exist.
\begin{corollary}[Global convergence]
  \label{co:globalconv}
  Let $\seq{(x\k,y\k,z\k)}$ be the sequence of vectors generated by
  Algorithm~\aref{algo:sLCL}.  Let $\xstar$ be the single limit point
  of the sequence $\seq \xstark$.  Suppose that
  Assumptions~\aref{ass:continuity}, \aref{ass:compactness},
  and~\aref{ass:licq} hold.  Then
  \begin{equation*}
    \lim_{k\to\infty} (x\k,y\k,z\k) = (\xstar,\ystar,\zstar),
  \end{equation*}
  and $(\xstar,\ystar,\zstar)$ is a first-order KKT point for \GNP.
\end{corollary}
\begin{proof}
  Let $\seq{(\xstark,\ystark,\zstark)}$ be the sequence of vectors
  generated by Step~\ref{algo:sLCLsolveELC} of
  Algorithm~\aref{algo:sLCL} and set $\fomult\k =
  \fomult(\xstark,\ystark,\rho\k)$.  By Lemma~\ref{le:penalty} and
  Theorem~\ref{th:globalconv},
  \begin{equation*}
    \lim_{k\to\infty}\fomult\k = \ystar
    \text{and}
    \lim_{k\to\infty}\zstark = \zstar.
  \end{equation*}
  Moreover, $(\xstar,\ystar,\zstar)$ is a first-order KKT point for
  \GNP.  Suppose Step~\ref{algo:sLCLupdatesolutions} is executed
  infinitely often.  The result then follows immediately because
  $x\k$, $y\k$, and $z\k$ are updated infinitely often and form a
  convergent sequence from $\xstark$, $\fomult\k$, and $\zstark$.
  
  We now show by contradiction that
  Step~\ref{algo:sLCLupdatesolutions} does occur infinitely often.
  Suppose instead that it does not.  Then there exists a $k_1$ large
  enough so that Steps~\ref{algo:sLCLkeep} and~\ref{algo:sLCLincrho}
  are executed for all $k > k_1$. Consider only iterations $k>k_1$.
  Then $y\k\equiv\ybar$ and $\rho\k\to\infty$.  From~\eqref{eq:70},
  \begin{equation}
    \label{eq:25}
  \begin{aligned}
    \rho\k \norm{c(\xstark)} &=   \norm{\ystark-\fomult\k}
\\                           &=   \norm{(\ystark-\ybar)+ \ybar-\fomult\k}
\\                           &\le \norm{\ystark-\ybar} +
                                  \norm{\ybar} +
                                  \norm{\fomult\k}.
  \end{aligned}
  \end{equation}
  The vector $\ystark\equiv y\k+\dystark$ and each $\dystark$
  satisfies~\eqref{eq:59-6}.  Thus,
  $\norm{\ystark-\ybar}\le\sigma\k+\omega\k$.  Moreover,
  $\lim_{k\to\infty}\fomult\k=\ystar$ and $\ystar$ is bounded (see
  Assumption~\ref{ass:licq}).  Then, from~\eqref{eq:25}, there exists
  some constant $L > 0$, independent of $k$, such that
  \begin{equation}
    \label{eq:8}
    \rho\k \norm{c(\xstark)} \le L
  \end{equation}
  for all $k$.  But the test at Step~\ref{algo:sLCLbranchtest} fails
  at every iteration, so that
  \begin{equation}
    \label{eq:71}
    \eta\k < \norm{c(\xstark)}.
  \end{equation}
  Combining~\eqref{eq:8} and~\eqref{eq:71}, we find that
  \begin{equation}
    \label{eq:72}
    \rho\k\eta\k < \rho\k \norm{c(\xstark)} \le L.
  \end{equation}
  From Step~\ref{algo:sLCLetaa}, $\eta\kp1=\eta_0 /
  \rho^{\alpha}\kp1$, so
  \begin{equation}
    \label{eq:74}
    \rho\k\eta\k
    = \rho\k \frac{\eta_0}{\rho\k^{\alpha}}
    = \eta_0 \rho\k^{1-\alpha}.
  \end{equation}
  Substituting~\eqref{eq:74} into~\eqref{eq:72}, we find that
  $\eta_0 \rho\k^{1-\alpha} < L$
  for all $k$.  This is a contradiction under the hypothesis
  that $\alpha<1$ and $\rho\k\to\infty$.  Therefore,
  Step~\ref{algo:sLCLupdatesolutions} must occur infinitely often.
\end{proof}

The following result simply asserts that Algorithm~\ref{algo:sLCL}
will eventually exit when, as in practice, $\omega_*$ and $\eta_*$ are
positive.
\begin{corollary}[Finite Termination]
  Suppose that the convergence tolerances $\omega_*$ and $\eta_*$ are
  strictly positive.  Then, under the assumptions of
  Corollary~\aref{co:globalconv}, Algorithm~\aref{algo:sLCL} terminates
  after a finite number of iterations.
\end{corollary}
\begin{proof}
  Let $\seq{(\xstark,\ystark,\zstark)}$ and $\xstar$ be as defined in
  Theorem~\ref{th:globalconv}.  Set $\fomult\k =
  \fomult\k(\xstark,\ystark,\rho\k)$.  By that theorem,
  \begin{align*}
    \lim_{k\to\infty} \fomult\k &= \ystar
\\  \lim_{k\to\infty} \zstark   &= \zstar \defd \grad{x}\AugL(\xstar,\ystar,0),
  \end{align*}
  and $(\xstar,\ystar,\zstar)$ is a first-order KKT point for \GNP.
  Then, $(\xstar,\ystar,\zstar)$ must satisfy~\eqref{eq:KKTGNPb}.  By
  the continuity of $c$, $\lim_{k\to\infty} \norm{c(\xstark)} \to
  \cstar = 0$, and because $\eta_*>0$,
  \begin{equation*}
    \norm{c(\xstark)} < \eta_* \le \max(\eta\k,\eta_*)
  \end{equation*}
  for all $k\in\Kscr$ large enough.  Consequently,
  Step~\ref{algo:sLCLtestconv} is executed infinitely often and
  \begin{equation*}
    \lim_{k\to\infty} (x\k,y\k,z\k) \to (\xstar,\ystar,\zstar).
  \end{equation*}
  Because $\omega_*>0$ and $\eta_*>0$, $(x\k,y\k,z\k)$ satisfies
  conditions~\eqref{eq:KKTGNPb} for some $k$ large enough.
\end{proof}

\subsection{Local convergence properties}        \label{sec:sLCL-local}

In this section we show that the stabilized LCL algorithm preserves
the local convergence characteristics of Robinson's original LCL
algorithm.  Moreover, it can retain fast local convergence under
inexact solutions to the subproblems.

Bertsekas~\cite{Ber82}\ and Conn~\etal~\cite{CGST96,CGT91b}\ show how
to construct a forcing sequence $\seq{\eta\k}$ to guarantee that
$\norm{c(\xstark)} \le \eta\k$ will eventually always be true so that
the iterates $x\k$, $y\k$, and $z\k$ are updated (see
Step~\ref{algo:sLCLupdatesolutions} of Algorithm~\ref{algo:sLCL}) for
all iterations after some $k$ large enough.  The penalty parameter
$\rho\k$ then remains uniformly bounded---an important property.
These results rely on a relationship between $\norm{c(\xstark)}$ and
$\rho\k$, namely~\eqref{eq:cconv1}.  We know from the BCL convergence
theory that the convergence rate approaches superlinear as $\rho\k$
grows large (\cf~\cite{Ber82} and~\cite{CGST96,CGT91b}).  Because
$\eta\k$ is reduced at a sublinear rate, $\norm{c(\xstark)}$ will
eventually go to zero faster than $\eta\k$, at which point it is no
longer necessary to increase $\rho\k$.  Thus, we can be assured that
Algorithm~\ref{algo:sLCL} does not increase $\rho\k$ without bound.

Bertsekas suggests constructing the sequence $\eta\k$ as
\begin{equation}
  \label{eq:75}
  \eta\kp1 =
    \gamma \norm{c(\xstark)},
\end{equation}
for some $\gamma < 1$.  Within Algorithm~\ref{algo:sLCL}, this would
lead to the following update rule:
\begin{equation}
  \label{eq:76}
  \rho\kp1 = 
  \begin{cases}
    \rho\k & \text{if $\norm{c(\xstark)} \le \gamma
      \norm{c(x\k)}$} \\
    \tau_\rho \rho\k & \text{if $\norm{c(\xstark)} > \gamma
      \norm{c(x\k)}$.}
  \end{cases}
\end{equation}
As $\rho\k$ gets larger, the convergence rate gets arbitrarily close
to superlinear, so that the first case of~\eqref{eq:76} is always
satisfied, and $\rho\k$ becomes constant for all $k$ large enough.  We
prefer not to use rule~\eqref{eq:75} because it may be too strict.
Any intermediate (and nonoptimal) iterate $\xstark$ could be feasible
or nearly feasible for \GNP, so that $\norm{c(\xstark)}$ could be very
small.  Then $\eta\kp1$ would be smaller than warranted on the
following iteration.  The forcing sequence suggested by Conn
\etal~\cite{CGST96,CGT91b}\ does not suffer from this defect and has
been proven by them to keep $\rho\k$ bounded.  We have used this
update in Algorithm~\ref{algo:sLCL} (see Steps~\ref{algo:sLCLetab}
and~\ref{algo:sLCLetaa}).

For this analysis and the remainder of this section, we assume that
$\rho\k$ is uniformly bounded, so that $\rho\k = \rhobar$ for all $k$
greater than some $\kbar$.  Hence, we drop the subscript on $\rho\k$
and simply write $\rhobar$.  We consider only iterations $k>\kbar$.

We begin by discussing the local convergence rates of the
Algorithm~\ref{algo:sLCL} under the assumption that the elastic
variables are always zero---that is, the linearized constraints are
always satisfied.  Next, we show that after finitely many iterations
the elastic penalty parameter $\sigma\k$ will always be large enough
to guarantee that this assumption holds.  In this way, we demonstrate
that stabilized LCL becomes equivalent to \MINOS\ (and to canonical
LCL) as it approaches the solution.

\subsubsection{Convergence rates}               \label{sec:sLCL-rates}

Robinson's~\cite{Rob72}\ local convergence analysis applies to the
canonical LCL algorithm under the special case in which $\rho\k \equiv
0$ (\cf \eqref{eq:4}) and each subproblem is solved to full accuracy
(\ie, $\omega\k \equiv 0$).  He proved that one can expect fast
convergence from a good enough starting point.  In particular, under
Assumptions~\ref{ass:continuity}, \ref{ass:second-order}, and
\ref{ass:licq}, we can expect an R-quadratic rate of convergence (see
Ortega and Rheinboldt~\cite{OR70}\ for an in-depth discussion of
root-convergence rates).  For a sufficiently good starting point,
Robinson~\cite{Rob74}\ proves that the subproblems \LCk\ are always
feasible.  He also shows that near a solution, the solutions to the LC
subproblems, if parameterized appropriately, form a continuous path
converging to $(\xstar,\ystar,\zstar)$.

In a later paper, Br\"{a}uninger~\cite{Bra77}\ shows how the fast
local convergence rate can be preserved with only \emph{approximate}
solutions of the subproblems (again, with $\rho\k\equiv0$).  The
subproblems are solved to a tolerance that is tightened at a rate that
matches the decrease in the square of the primal and dual
infeasibilities.  Our proposed LCL algorithm uses a similar strategy.

Robinson's local convergence analysis also applies to the canonical
LCL algorithm when $\rho\k\equiv\rhobar>0$.  One can see this by
considering the following optimization problem:
\begin{equation}
  \label{eq:77}
  \begin{array}{ll}
    \minimize{x} & f(x) + \half \rhobar \norm{c(x)}^2 \\[6pt]
    \st          & c(x) = 0, \quad x \ge 0.
  \end{array}
\end{equation}
The solutions of~\eqref{eq:77} are identical to the solutions of \GNP.
The Robinson LCL subproblem objective corresponding to~\eqref{eq:77}
is given by
\begin{equation*}
  R\k(x) = f(x) + \half \rhobar \norm{c(x)}^2 - y\k\T c(x).
\end{equation*}
The canonical LCL subproblem objective is $\AugL\k(x) \equiv
\AugL(x,y\k,\rho\k)$, and so $\AugL\k(x) \equiv R\k(x)$ for all $k$
because $\rho\k \equiv \rhobar$.  We then observe that the canonical
LCL subproblem corresponding to \GNP, with a penalty parameter $\rho\k
\equiv \rhobar$, is \emph{equivalent} to the Robinson LCL subproblem
corresponding to problem~\eqref{eq:77}, with $\rho\k\equiv0$.  The
convergence characteristics of the canonical LCL algorithm are
therefore the same as those demonstrated by Robinson~\cite{Rob72}.
(However, while the asymptotic convergence rate remains R-quadratic,
we expect a different asymptotic error constant.)

Under the assumption that the elastic variables are always equal to 0
and that $\rhobar$ is finite, the steps executed by
Algorithms~\ref{algo:cLCL} and~\ref{algo:sLCL} are identical, and the
subproblems \ELCk\ and \LCk\ are also identical.  The only difference
is the multiplier update formulas:
\begin{subequations}
\begin{align}
  \label{eq:78}
  \text{Canonical  LCL update} &\quad y\kp1 = \ystark \\
  \text{Stabilized LCL update} &\quad y\kp1 = \ystark - \rhobar
  c(\xstark),
\end{align}
\end{subequations}
which differ only by the vector $\rhobar c(\xstark)$.  We may think of
this vector as a perturbation of the LCL multiplier
update~\eqref{eq:78}.  Moreover, Robinson~\cite{Rob72} shows that this
perturbation converges to 0 at the same rate as $\seq \xstark$
converges to $\xstar$.  Therefore, it does not interfere with the
convergence rate of the stabilized LCL iterates.  Robinson's local
convergence analysis then applies to the stabilized LCL method.

We summarize the convergence results in Theorem~\ref{th:lclrates}.
Note that the function
\begin{equation*}
  F(x,y,z) =
  \begin{bmatrix}
       c(x)
    \\ \grad{x}\AugL(x,y,\rho) - z
    \\ \min(x,z)
  \end{bmatrix}
 \end{equation*}
 captures the first-order optimality conditions of \GNP, in the sense
 that $$F(\xstar,\ystar,\zstar) = 0$$
 if and only if
 $(\xstar,\ystar,\zstar)$ is a first-order KKT point for \GNP.  Thus,
 $\norm{F(x,y,z)}$ is a measure of the deviation from optimality.  For
 the next theorem only, define
 \begin{equation*}
   r =
   \begin{pmatrix}
     x \\ y \\ z
   \end{pmatrix}
   \text{and}
   F(r) = F(x,y,z).
 \end{equation*}
\begin{theorem}[Robinson~\cite{Rob72}; Br\"{a}uninger~\cite{Bra77}]
  \label{th:lclrates}
  Suppose Assumptions~\aref{ass:continuity}--\aref{ass:second-order}
  and~\aref{ass:licq} hold at $\rstar$.  Moreover, suppose $\omega\k =
  O(\norm{F(r\k)}^2)$ for all $k \ge 0$.  Then there is a positive
  constant $\delta$ such that if
  \begin{equation*}
    \norm{r_0 - r_*} < \delta,
  \end{equation*}
  the sequence $\seq{r\k}$ generated by Algorithm~\aref{algo:sLCL}
  converges to $\rstar$.  Moreover, the sequence converges
  $R$-quadratically, so that for all $k \ge 0$,
\begin{equation}
  \label{eq:80}
  \norm{r\k - \rstar} \le Q (\half)^{2^k}
\end{equation}
for some positive constant $Q$.  Also,
\begin{equation}
  \label{eq:79}
  \norm{r\kp1 - r\k} \le M \norm{F(r\k)}
\end{equation}
for some positive constant $M$.
\end{theorem}

Robinson does not state~\eqref{eq:79} as part of a theorem, but it is
found in the proof of~\eqref{eq:80}.

\subsubsection{Asymptotic equivalence to MINOS}  \label{sec:sLCL-equiv}

Much of the efficiency of LCL methods, including \MINOS, derives from
the fact that they eventually identify the correct active set, and
each subproblem restricts its search to the subspace defined by a
linear approximation of the constraints.  This approximation can be
very accurate near the solution.  The stabilized LCL subproblems do
\emph{not} restrict themselves to this subspace.  In early iterations
we do not expect, nor do we wish, the method to honor these
linearizations.  The elastic variables give the subproblems an
opportunity to deviate from this subspace.  In order to recover LCL's
fast convergence rate, however, it is not desirable to allow deviation
near the solution.

We show below that as the stabilized LCL iterations approach a
solution of \GNP, the solutions of the stabilized LCL subproblems
eventually always have the elastic variables equal to zero.  Hence,
$\clink{\xstark} = 0$ and $\vstark,\wstark=0$, so that each $\xstark$
satisfies the constraints of the canonical LCL (and \MINOS)
subproblem, and the objective of \ELCk\ at $(\xstark,\vstark,\wstark)$
is equivalent to \LCk.

As discussed in \Sec\ref{sec:sLCL-rates}, Theorem~\ref{th:lclrates}
applies to Algorithm~\ref{algo:sLCL}, and so for all $k$ large enough,
\eqref{eq:79} yields
\begin{equation}
  \label{eq:35}
  \norm{\ystark - y\k} \le M \norm{F(x\k,y\k,z\k)},
\end{equation}
where $\ystark\equiv y\k+\dystark$ and $M$ is some positive constant.
By Corollary~\ref{co:globalconv}, $(x\k,y\k,z\k) \to
(\xstar,\ystar,\zstar)$, and because $\norm{F(\xstar,\ystar,\zstar)} =
0$ and $F$ is continuous,
\begin{equation}
  \label{eq:36}
  \norm{F(x\k,y\k,z\k)} < \frac{\sigmab}{M}
\end{equation}
for all $k$ large enough ($\sigmab$ is defined in
Algorithm~\ref{algo:sLCL}).  Combining~\eqref{eq:35}
and~\eqref{eq:36}, we conclude that $\norm{\ystark - y\k} < \sigmab$,
or equivalently,
\begin{equation}
  \label{eq:38}
  \norm{\dystark} < \sigmab
\end{equation}
for all $k$ large enough.  However, Step~\ref{algo:sLCLupdatesigma} of
Algorithm~\ref{algo:sLCL} guarantees that $\sigmab \le \sigma\k$ for
all $k$, and so from~\eqref{eq:38}, $\norm{\dystark} < \sigma\k$ for
all $k$ large enough.  Lemma~\ref{le:L1threshold} then implies that
$\sigma\k$ will be sufficiently large that the optimal elastic
variables will be equal to 0.

\subsection{Infeasible problems}                \label{sec:sLCL-infeas}

Not all optimization problems are feasible.  The user of an
optimization algorithm may formulate a set of nonlinear constraints
$c(x)=0$ for which no nonnegative solution exists.  Detecting
infeasibility of the system $c(x) = 0$, $x \ge 0$, is equivalent to
verifying that the \emph{global} minimizer of
\begin{equation}
  \label{eq:56}
  \begin{array}{ll}
    \minimize{x} & \half \norm{c(x)}^2 \\
    \st & x \ge 0
  \end{array}
\end{equation}
yields a positive objective value.  Detecting such infeasibility
is a useful feature, but it is a very difficult
problem and is beyond the purview of this paper.

We analyze the properties of the stabilized LCL algorithm when it is
applied to an infeasible problem with convergence tolerances
$\omega_*=\eta_* = 0$.  We show that Algorithm~\ref{algo:sLCL}
converges to a point that satisfies the first-order optimality
conditions of the minimum-norm problem~\eqref{eq:56}.
\begin{theorem}
  \label{th:GNPinfeas}
  Let $\xstar$ be any limit point of the sequence of vectors $\seq
  \xstark$ generated by Algorithm~\aref{algo:sLCL}, and let $\Kscr$ be
  the infinite set of indices associated with that subsequence.
  Suppose that \GNP\ is infeasible.  Then, under the assumptions of
  Lemma~\aref{le:penalty},
  \begin{equation*}
   \lim_{k\in\Kscr}J(\xstark)\T c(\xstark) = \zstar \defd \Jstar\T \cstar,
  \end{equation*}
  and $(\xstar,\zstar)$ is a first-order KKT point for~\eqref{eq:56}.
\end{theorem}
\begin{proof}
  The pair $(\xstar,\zstar)$ satisfies the first-order KKT conditions
  of~\eqref{eq:56} if
  \begin{equation}
  \label{eq:14}
  \begin{aligned}
    \Jstar\T \cstar       &=   \zstar
\\  \min(\xstar,\zstar)   &=   0.
  \end{aligned}
  \end{equation}
  Because \GNP\ is infeasible, there exists a constant $\delta > 0$
  such that $\delta < \norm{c(x)}$ for all $x\ge 0$.  Moreover,
  Steps~\ref{algo:sLCLetab} and~\ref{algo:sLCLetaa} of
  Algorithm~\ref{algo:sLCL} generate a sequence $\seq{\eta\k}$
  converging to 0, and so $\eta\k < \delta$ for all $k$ large enough.
  Consider only such $k$.  Then, $\eta\k < \delta <
  \norm{c(\xstark)}$, and Step~\ref{algo:sLCLincrho} is executed at
  every $k$, so that $\rho\k\to\infty$ and $\sigma\k \to 0$.
  Moreover, $x\k$ and $y\k$ are not updated, so that for some
  $n$-vector $\xbar$ and $m$-vector $\ybar$,
  \begin{equation}
    \label{eq:61}
    x\k \equiv \xbar \text{and} y\k \equiv \ybar.
  \end{equation}
  
  Note that Algorithm~\ref{algo:sLCL} generates $\xstark$
  satisfying~\eqref{eq:59}.  Therefore, $\xstark \ge 0$ for all $k$,
  and so $\lim_{k\in\Kscr}\xstark = \xstar$ implies
  \begin{equation}
    \label{eq:11}
    \xstar\ge0.
  \end{equation}

  From~\eqref{eq:59-d}, \eqref{eq:40}, and~\eqref{eq:61},
  \begin{equation}
    \label{eq:60}
    g(\xstark) - J(\xstark)\T (\ybar - \rho\k c(\xstark))
    - J(\xbar)\T \dystark \ge -\omega\k e,
  \end{equation}
  or, after rearranging terms,
  \begin{equation}
    \label{eq:62}
    \underbrace{g(\xstark) - J(\xstark)\T\ybar}_{(a)}
    - \underbrace{J(\xbar)\T \dystark}_{(b)}
    + \rho\k J(\xstark)\T c(\xstark) \ge -\omega\k e.
  \end{equation}
  By hypothesis, all iterates $\xstark$ lie in a compact set, and so
  $(a)$ is bounded because $g$ and $J$ are continuous and $\ybar$ is
  constant.  Also, $(b)$ is bounded because $\xbar$ is constant, and
  from~\eqref{eq:59-6} we have $\infnorm{\dystark} \le \sigma\k +
  \omega\k$.  Then, because $\omega\k \to 0$ and $\rho\k \to \infty$,
  \eqref{eq:62} implies that $J(\xstark)\T c(\xstark) \ge 0$ for all
  $k$ large enough.  Otherwise,~\eqref{eq:62} would eventually be
  violated as $\rho\k$ grew large.  Then,
  \begin{equation}
    \label{eq:12}
    \zstar\defd\lim_{k\in\Kscr} J(\xstark)\T
    c(\xstark) = \Jstar\T \cstar \ge 0.
  \end{equation}
  
  All $\xstark$ lie in a compact set, there exists some constant $L>0$
  such that
  \begin{equation}
    \label{eq:81}
    \norm{\xstark-\xbar} \le \frac{L\alpha_1}{\sqrt n M},
  \end{equation}
  where $M$ and $\alpha_1$ are as defined in Lemma~\ref{le:penalty}
  and $n$ is the number of elements in the vector $\xstark$.
  Substituting~\eqref{eq:81} into~\eqref{eq:67} and
  using~\eqref{eq:23}, we have
  \begin{equation}
    \label{eq:73}
    \norm{\ghat(\xstark) - \Jhat(\xstark)\T\ystark
      + \rho\k \Jhat(\xstark)\T c(\xstark)}
    \le
    \sqrt n \{\omega\k + L(\sigma\k+\omega\k)\}.
  \end{equation}
  Dividing~\eqref{eq:73} through by $\rho\k$, we obtain
  \begin{equation}
    \label{eq:83}
    \left\|\frac{1}{\rho\k}
      \bigl(\ghat(\xstark) - \Jhat(\xstark)\T\ystark\bigr)
      + \Jhat(\xstark)\T c(\xstark)\right\|
    \le
    \frac{\sqrt n\{\omega\k + L(\sigma\k+\omega\k)\}}{\rho\k}.
  \end{equation}
  The quantity $\ghat(\xstark) - \Jhat(\xstark)\T\ystark$ is bounded
  for the same reasons that $(a)$ and $(b)$ above are bounded.  Taking
  limits of both sides of~\eqref{eq:83}, $\rho\k\to\infty$ and
  $\omega\k,\sigma\k\to 0$ imply that $\Jhat(\xstark)\T c(\xstark) \to
  0$.  By continuity of $J$ and $c$, $\Jhat_*^{\; T} \cstar = 0$.
  Equivalently, we may write
  \begin{equation}
    \label{eq:13}
    \compj{\Jstar\T \cstar} = 0 \text{if} \compj{\xstar} > 0,
  \end{equation}
  for $j = 1,\ldots,n$.  Therefore \eqref{eq:11}, \eqref{eq:12}
  and~\eqref{eq:13} together imply that $(\xstar,\zstar)$ satisfies
  conditions~\eqref{eq:14}, as required.
\end{proof}

Theorem~\ref{th:GNPinfeas} describes a useful feature of
Algorithm~\ref{algo:sLCL}.  When applied to an infeasible problem, the
algorithm converges to a solution of~\eqref{eq:56}---or at least to a
first-order point.  One important caveat deserves mention: if the
convergence tolerance $\eta_*$ is small (it usually will be),
Algorithm~\ref{algo:sLCL} may never terminate.  We need to insert an
additional test to provide for the possibility that \GNP\ is
infeasible.  For example, the test could force the algorithm to exit
if $\rho\k$ is above a certain threshold value and $\norm{c(\xstark)}$
is no longer decreasing.  Any test we devise is necessarily heuristic,
however; it is impossible to know for certain whether a larger value
of $\rho\k$ would force $\norm{c(\xstark)}$ to be less than $\eta_*$.
We discuss this point further in \Sec\ref{sec:implem-infeas}.

\subsection{Second-order optimality}            \label{sec:sLCL-second}

The stabilized LCL method imposes few requirements on the manner in
which the LC subproblems are solved.  Our implementation (see
Section~\ref{sec:imp}) uses \MINOS\ or \SNOPT\ to solve the LC
subproblems.  These are active-set solvers suitable for optimization
problems with few expected degrees of freedom at the solution and in
which only first derivatives are available.  However, second
derivatives might be readily available for some problems.  Also, some
problems are expected to have many degrees of freedom at the solution.
In either case, an interior-point solver (requiring second
derivatives) may be more appropriate for the solution of the
subproblems.

Lemma~\ref{le:penalty} and Theorem~\ref{th:globalconv} assert that
iterates generated by the stabilized LCL algorithm converge to
first-order KKT points.  A subproblem solver that uses
second-derivatives may be able to guarantee convergence to
second-order points.  If we augment the convergence criteria for the
solution of each subproblem to include second-order conditions, we can
show that Algorithm~\ref{algo:sLCL} generates iterates converging to
points satisfying the second-order sufficiency conditions for \GNP.
The following assumption strengthens the first-order
conditions~\eqref{eq:58}.
\begin{assumption}
  \label{ass:ELCsecond-order}
  Let $\xstar$ be any limit point of the sequence $\seq\xstark$, and
  let $\Kscr$ be the infinite set of indices associated with that
  convergent subsequence.  For all $k\in\Kscr$ large enough, the
  following conditions hold at each $(\xstark,\ystark,\zstark)$: For
  some $\delta > 0$, independent of $k$,
  \begin{enumerate}
  \item (Strict Complementarity)
    \begin{equation}
      \label{eq:91}
      \max(\xstark,\zstark) > \delta e;
    \end{equation}
  \item (Second-Order Condition) For any
    $\rho\ge0$,
    \begin{equation}
      \label{eq:89}
      p\T \Hess{xx}\AugL(\xstark,\ystark,\rho)p > \delta
    \end{equation}
    for all $p\neq0$ satisfying
    \begin{equation}
      \label{eq:15}
      \text{$J(\xstark) p=0$ and $\compj p=0$ for all $j$ such that
      $\compj \xstark=0$.}
  \end{equation}
  \end{enumerate}
\end{assumption}
\noindent
Condition~\eqref{eq:89} implies that the reduced Hessian of $\AugL$
is uniformly positive definite at all $\xstark$.

The following result extends Theorem~\ref{th:globalconv} to consider
the case in which iterates generated by Algorithm~\ref{algo:sLCL}
satisfy Assumption~\ref{ass:ELCsecond-order}.  Conn
\etal~\cite{CGT91b}\ show a similar result for their BCL method.
\begin{theorem}
  \label{th:second-order}
  Suppose that Assumptions~\aref{ass:continuity},
  \aref{ass:compactness}, \aref{ass:licq},
  and~\aref{ass:ELCsecond-order} hold.  Let
  $\seq{(\xstark,\ystark,\zstark)}$ be the sequences of vectors
  generated by Algorithm~\aref{algo:sLCL}.  Let $\xstar$ be any limit
  point of the sequence $\seq \xstark$, and let $\Kscr$ be the
  infinite set of indices associated with that convergent subsequence.
  Set $\fomult\k = \fomult(\xstark,\ystark,\rho\k)$.  Then
  \begin{equation}
    \label{eq:92}
    \lim_{k\in\Kscr} (\xstark,\fomult\k,\zstark) = (\xstar,\ystar,\zstar)
  \end{equation}
  and $(\xstar,\ystar,\zstar)$ is an isolated local minimizer of \GNP.
\end{theorem}
\begin{proof}
  It follows immediately from Theorem~\ref{th:globalconv}
  that~\eqref{eq:92} holds and that $(\xstar,\ystar,\zstar)$ is a
  first-order KKT point for \GNP.  It only remains to show that
  $(\xstar,\ystar,\zstar)$ satisfies the second-order sufficiency
  conditions (see Definition~\ref{def:second-order}).
  
  By hypothesis, $\xstark$ and $\zstark$ satisfy Part~1 of
  Assumption~\ref{ass:ELCsecond-order} for all $k\in\Kscr$.
  Therefore, their limit points satisfy
  \begin{equation*}
    \max(\xstar,\zstar) \ge \delta e > 0,
  \end{equation*}
  and so $\xstar$ and $\zstar$ satisfy strict complementarity
  (Definition~\ref{def:strict-comp}).  We now show that $\xstar$ and
  $\ystar$ satisfy the second-order sufficiency conditions for \GNP.
  
  Let $p$ be any nonzero vector satisfying~\eqref{eq:15} for all
  $k\in\Kscr$ large enough.  Then
  \begin{equation}
    \label{eq:17}
    \begin{aligned}
      p\T\Hess{xx}\AugL(\xstark,\ystark,\rho\k) p
      &= p\T\big(
         H(\xstark) - \sum_{i=1}^{m}\compi{\fomult\k}H_i(\xstark)
         \big)p
    \end{aligned}
  \end{equation}
  for all $k\in\Kscr$ large enough.  Part~2 of
  Assumption~\ref{ass:ELCsecond-order} and~\eqref{eq:17} imply that
  \begin{equation}
    \label{eq:19}
    p\T\big(
         H(\xstark) - \sum_{i=1}^{m}\compi{\fomult\k}H_i(\xstark)
         \big)p > \delta,
  \end{equation}
  where $\delta$ is some positive constant.  If we take the limit
  of~\eqref{eq:19}, the continuity of $H$ and $H_i$ (see
  Assumption~\ref{ass:continuity}) and \eqref{eq:92} imply that
  \begin{equation}
    \label{eq:24}
    p\T \Hess{xx}\AugL(\xstar,\ystar,\rho)p =
    p\T\big(
         H(\xstar) - \sum_{i=1}^{m}\compi\ystar H_i(\xstar)
    \big)p
    \ge \delta > 0
  \end{equation}
  for all $\rho\ge0$ and for all $p\ne0$ satisfying~\eqref{eq:28}.
  Therefore, $(\xstar,\ystar,\zstar)$ satisfies the second-order
  sufficiency conditions for \GNP, as required.
\end{proof}

\section{Implementation}                                \label{sec:imp}

The practical implementation of an algorithm invariably requires many
features that are not made explicit by its theory.  In this section we
discuss some important details of our implementation of the stabilized
LCL method.  The algorithm has been implemented in \Matlab,
version~6~\cite{Matlabuser}\ and is called \LCLOPT.  It uses the
Fortran codes \MINOS~\cite{MS78,MS82}\ and \SNOPT~\cite{GMS97b}\ to
solve the linearly constrained subproblems.  We now turn our attention
back to the more general problem \NP, first presented in
\Sec\ref{sec:intro-nlp}, and leave \GNP\ behind.

\subsection{Problem formulation}                   \label{sec:implem-form}

\LCLOPT\ does not solve \NP\ directly, but rather solves the
equivalent problem
\begin{equation*}
  \problem{\NPi}{x,s}
  {f(x)}
  {
    \begin{pmatrix}
      c(x) \\ Ax
    \end{pmatrix}
    - s
    = 0,
    \quad
    l \le
    \begin{pmatrix}
      x \\ s
    \end{pmatrix}
    \le u.
    \vspace{6pt}
  }
\end{equation*}
The formulation of problem \NPi\ is chosen to match the problem
formulation used by \SNOPT.  It is also closely related to that used
by \MINOS.  As in those methods, our implementation distinguishes
between variables in the vector $x$ that appear and do not appear
nonlinearly in the objective or the constraints; variables that appear
only linearly are treated specially.  The following discussion ignores
this detail in order to keep the notation concise.

The linearly constrained subproblems corresponding to \NPi\ take the
form
\begin{equation*}
  \wideproblem{\ELCik}{x,s,v,w}
  {\AugL\k(x) + \sigma\k e\T (v+w)}
  {\begin{aligned}[t]
      \begin{pmatrix}
        c\k + J\k(x - x\k) + v - w \\ Ax
      \end{pmatrix}
      - s = 0,
      \quad
    l &\le
    \begin{pmatrix}
      x \\ s
    \end{pmatrix}
    \le u,
\\  0 &\le v,w.
    \end{aligned}
  }
\end{equation*}

\subsection{The main algorithm}                 \label{sec:implem-main}

The computational kernel of \LCLOPT\ resides in the solution of each
LC subproblem, and the efficiency of the implementation ultimately
relies on the efficiency of the subproblem solver.  The main tasks of
the outer level are to form the subproblems, update solution
estimates, update parameters, and test for convergence or errors.

\subsection{Solving the LC subproblems}      \label{sec:implem-LCprobs}

\LCLOPT\ can use either \MINOS\ or \SNOPT\ to solve \ELCik.  For
linearly constrained problems, \MINOS\ uses a reduced-gradient method,
coupled with a quasi-Newton approximation of the reduced Hessian of
the the problem objective.  \SNOPT\ implements a sparse SQP method and
maintains a limited-memory, quasi-Newton approximation of the Hessian
of the problem objective. (In both cases, the problem objective will
be the objective of \ELCik.)  For linearly constrained problems,
\SNOPT\ avoids performing an expensive Cholesky factorization of the
reduced Hessian for the quadratic programming subproblem in each of
its own major iterations, and thus realizes considerable computational
savings over problems with nonlinear constraints~\cite{GMS02}.

Both \MINOS\ and \SNOPT\ are available as libraries of Fortran 77
routines.  We implemented MEX interfaces~\cite{Matlabapi} written in C
to make each of the routines from the \MINOS\ and \SNOPT\ libraries
accessible from within \Matlab.  The subproblem solvers evaluate the
nonlinear objective function (there are no nonlinear constraints in
\ELCik) through a generic MEX interface, \texttt{funObj.c}.  This
routine makes calls to a \Matlab\ routine to evaluate the nonlinear
objective $\AugL\k$.  In turn, the routine for $\AugL\k$ makes calls
to routines (available as \Matlab\ or MEX routines) to evaluate the
original nonlinear functions $f$ and $c$.

\subsection{Computing an initial point}         \label{sec:implem-init}

\MINOS\ and \SNOPT\ both ensure that all iterates remain feasible (to
within a small tolerance) with respect to the bounds and linear
constraints in \ELCik, which includes the bounds and linear
constraints in \NPi.  \LCLOPT\ is therefore able to restrict the
evaluation of the nonlinear functions $f$ and $c$ to points in the
latter region.  A user of \LCLOPT\ may thus introduce bounds and
linear constraints into \NPi\ to help guard against evaluation of the
nonlinear functions at points where they are not defined.

Before entering the first iteration of the stabilized LCL method,
\LCLOPT\ solves the following quadratic \emph{proximal-point} (PP)
problem:
\begin{equation*}
  \problem{\PP2}{x}{\half\twonorm{x - \xtilde}^2}
  {l \le
    \begin{pmatrix}
      x \\ Ax
    \end{pmatrix}
   \le u,\vspace{6pt}}
\end{equation*}
where $\xtilde$ is a vector provided by the \LCLOPT\ user.  The
solution of \PP2\ is used as the initial point $x_0$ for the
algorithm.  The objective function of the PP problem helps find an
$x_0$ reasonably close to $\xtilde$, while the constraints ensure that
$x_0$ is feasible with respect to the bounds and linear constraints of
\NPi.  If \PP2\ proves infeasible, \NPi\ is declared infeasible and
\LCLOPT\ exits immediately with an error message.

An alternative PP problem is based on the one-norm deviation from
$\xtilde$:
\begin{equation*}
  \problem{\PP1}{x}{\onenorm{x - \xtilde}}
  {l \le
    \begin{pmatrix}
      x \\ Ax
    \end{pmatrix}
   \le u.\vspace{6pt}}
\end{equation*}
An advantage of \PP1\ is that it can be reformulated and solved as a
linear program, and its solution is therefore expected to lie on more
constraint vertices.  It is a correspondingly easier problem to solve
for reduced-space solvers.  \SNOPT\ provides the option of solving
either \PP1\ or \PP2.  \LCLOPT\ can take advantage of this by
initializing $x=\xtilde$ and passing \SNOPT\ the optimization problem
\begin{equation*}
  \begin{array}{ll}
    \minimize{x} & 0 \\
    \st & 
    l \le
    \begin{pmatrix}
      x \\ Ax
    \end{pmatrix}
   \le u,
  \end{array}
\end{equation*}
with the parameter \opt{Proximal Point} set to either 1 or 2.  (For
\MINOS, the constraints would have to be reformulated and a set of
elastic variables introduced.)

The computational results presented in \Sec\ref{sec:num} were derived
by using \PP2\ to compute $x_0$.  As suggested by Gill
\etal~\cite{GMS02}, a loose optimality tolerance on \PP2\ is used to
limit the computational expense of its solution: reducing the number
of iterations and (typically) the number of superbasic variables.

\subsection{Early termination of the LC subproblems}
\label{sec:implem-early-term}

The global convergence results for the stabilized LCL algorithm
(\cf~Lemma~\ref{le:penalty} and Theorem~\ref{th:globalconv}) assume
that the optimality tolerances $\omega\k$ for the subproblems
converge to 0.  This loose requirement allows much
flexibility in constructing the sequence $\seq{\omega\k}$.

The solution estimates may be quite poor during early iterations.  We
expect slow progress during those iterations, even if they are solved
to tight optimality tolerances.  A loose tolerance may help limit the
computational work performed by the subproblem solver during these
early iterations.  Near a solution, however, we wish to reduce the
optimality tolerance quickly in order to take advantage of the fast
local convergence rate predicted by Theorem~\ref{th:lclrates}.

To construct the sequence \seq{\omega\k}, we replace
Step~\ref{algo:sLCLomega} of Algorithm~\ref{algo:sLCL} by
\begin{equation}
  \label{eq:97}
  \begin{aligned}
    \omega     &\assign \min(\omega\k,\infnorm{F(x\k,y\k,z\k)}^2)
\\  \omega\kp1 &\assign \max(0.5 \omega,\omega_*),
  \end{aligned}
\end{equation}
where $\omega_0$ can be set by a user to any value between 0.5 and
$\omega_*$.  The update~\eqref{eq:97} guarantees that $\omega\k \to
\omega_*$, as required.

Following the prescription outlined in \Sec\ref{sec:sLCL-earlyterm},
we fix at a small value the feasibility tolerance for satisfying the
linearized constraints.  The feasibility and optimality tolerances for
each major iteration are passed to the subproblem solver as run-time
parameters.

\subsection{Detecting infeasibility and unboundedness}\label{sec:implem-infeas}

As discussed in \Sec\ref{sec:sLCL-infeas}, Algorithm~\ref{algo:sLCL}
will \emph{not} exit if the optimization problem is infeasible and the
infeasibility tolerance $\eta_*$ is small.  We declare \NPi\ 
infeasible if at any given iteration $k$, $x\k$ is infeasible with
respect to the nonlinear constraints and the penalty parameter is
greater than some threshold value.  In particular, at
Step~\ref{algo:sLCLincrho}, Algorithm~\ref{algo:sLCL} exits and \NPi\ 
is declared infeasible if
\begin{equation*}
  \begin{aligned}
    \max(\infnorm{\pospart{l_c - c\k}},
         \infnorm{\pospart{u_c - c\k}}) &> \eta_*
\\  \rho\k &> \rhobar,
  \end{aligned}
\end{equation*}
where $l_c$ and $u_c$ are the lower and upper bounds for the nonlinear
constraints and $\pospart{\cdot}$ is the positive part of a vector.
For the computational results in \Sec\ref{sec:num} the threshold value
was set at $\rhobar= 10^8$.

We also need to consider the possibility that \NPi\ is
unbounded---\ie, that the objective function $f$ is unbounded below in
the feasible region, or that $\norm x \to \infty$.  As with tests for
infeasibility, any test for unboundedness must be ad hoc.  We rely on
the LC solver to help detect infeasibility.  Problem \NPi\ is declared
unbounded and \LCLOPT\ exits if the point $x\k$ is feasible and the LC
solver reports \ELCik\ as unbounded.

\subsection{Summary of the stabilized LCL method}\label{sec:implem-summ}

Following is a summary of the stabilized LCL method as implemented in
\LCLOPT.  We assume that $\xtilde$ is given and that the starting
tolerances, $\omega_0$ and $\eta_0$, and parameters, $\rho_0$ and
$\sigma_0$, are set.
\begin{enumerate}
\item Apply the LC solver to \PP1\ or \PP2\ to obtain a starting point
  $x_0$ that is feasible with respect to the bounds and linear
  constraints and reasonably close to $\xtilde$.  If the PP problem is
  infeasible, declare \NPi\ infeasible and exit.  Otherwise, set
  $k=0$.
\item \label{step:linearize} Evaluate the functions and gradients at
  $x\k$.  Linearize the constraints and form \ELCik.
\item Apply the LC solver to \ELCik\ with optimality tolerance $\omega\k$
  to obtain $(\xstark,\dystark,\zstark)$.
  Set $\ystark = y\k + \dystark$.
\item If \ELCik\ is unbounded and $\xstark$ is \emph{feasible},
  declare \NPi\ unbounded and exit.  If \ELCik\ is unbounded and
  $\xstark$ is \emph{infeasible}, go to Step~\ref{step:fail}.
  Otherwise, continue.
\item If $\xstark$ meets the current nonlinear feasibility threshold
  $\eta\k$, continue.  Otherwise, go to Step~\ref{step:fail}.
\item \label{step:success} Update the solution estimate:
  $(x\kp1,y\kp1,z\kp1)\assign(\xstark,\ystark - \rho\k c(\xstark),\zstark)$.
  Keep the penalty parameter $\rho\k$ fixed and reset the elastic weight
  $\sigma\k$.
\item Test convergence: If $(x\kp1,y\kp1,z\kp1)$ satisfies the
  optimality conditions for \NPi, declare the current solution
  estimate optimal, return $(x\kp1,y\kp1,z\kp1)$, and exit.
  Otherwise, go to Step~\ref{step:tols}.
\item \label{step:fail} If $\rho\k>\rhobar$, declare \NPi\ infeasible,
  return $(\xstark,\ystark,\zstark)$, and exit.  Otherwise, discard
  the subproblem solution (\ie,
  $(x\kp1,y\kp1,z\kp1)\assign(x\k,y\k,z\k)$), increase the penalty
  parameter $\rho\k$, and reduce the elastic weight $\sigma\k$.
\item \label{step:tols} Set the next nonlinear feasibility threshold
  $\eta\kp1$ and LC subproblem optimality tolerance $\omega\kp1$, so
  that $\{(\omega\k,\eta\k)\}\to(\omega_*,\eta_*)$.
\item \label{step:iterate} Set $k \assign k+1$.  Return to
  Step~\ref{step:linearize}.
\end{enumerate}

\section{Numerical Results}                             \label{sec:num}

This section summarizes the results of applying our implementation of
the stabilized LCL method, \LCLOPT, to a subset of nonlinearly
constrained test problems from the \COPS~2.0~\cite{DM00},
Hock-Schittkowski~\cite{HS81}, and \CUTE~\cite{BCGT95}\ test suites.
Two versions of \LCLOPT\ are applied to each test problem: The first
version uses \AMPL/\MINOS~5.5~\cite{FGK93}, version 19981015, to solve
the sequence of linearly constrained subproblems; the second version
uses \SNOPT\ version~6.1-1(5).

We used the \AMPL\ versions of all problems, as formulated by
Vanderbei~\cite{Van02}.  A MEX interface to the \AMPL\ libraries makes
functions and gradients available in \Matlab\ (see Gay~\cite{Gay97}\ 
for details on interfacing external routines to \AMPL).  All runs were
conducted on an AMD Athlon 1700XP using 384 MB of RAM, running Linux
2.4.18.

\begin{figure}[htbp]
  \label{fig:perfprof}
  \centering
  \small
  \vspace*{-10pt}
  \psfrag{ylabel}[]{\% of problems}
  \psfrag{xlabel}[]{$\tau$ $\times$ best performer}
  \psfrag{title}[]{}
  \psfrag{tab-lclmi}[br][][.7]{\tiny \LCLOPT\ with \MINOS}
  \psfrag{tab-lclsn}[br][][.7]{\tiny \LCLOPT\ with \SNOPT}
  \psfrag{tab-minos}[br][][.7]{\tiny \MINOS}
  \subfigure[Function evaluations]%
  {\includegraphics[height=2.3in]{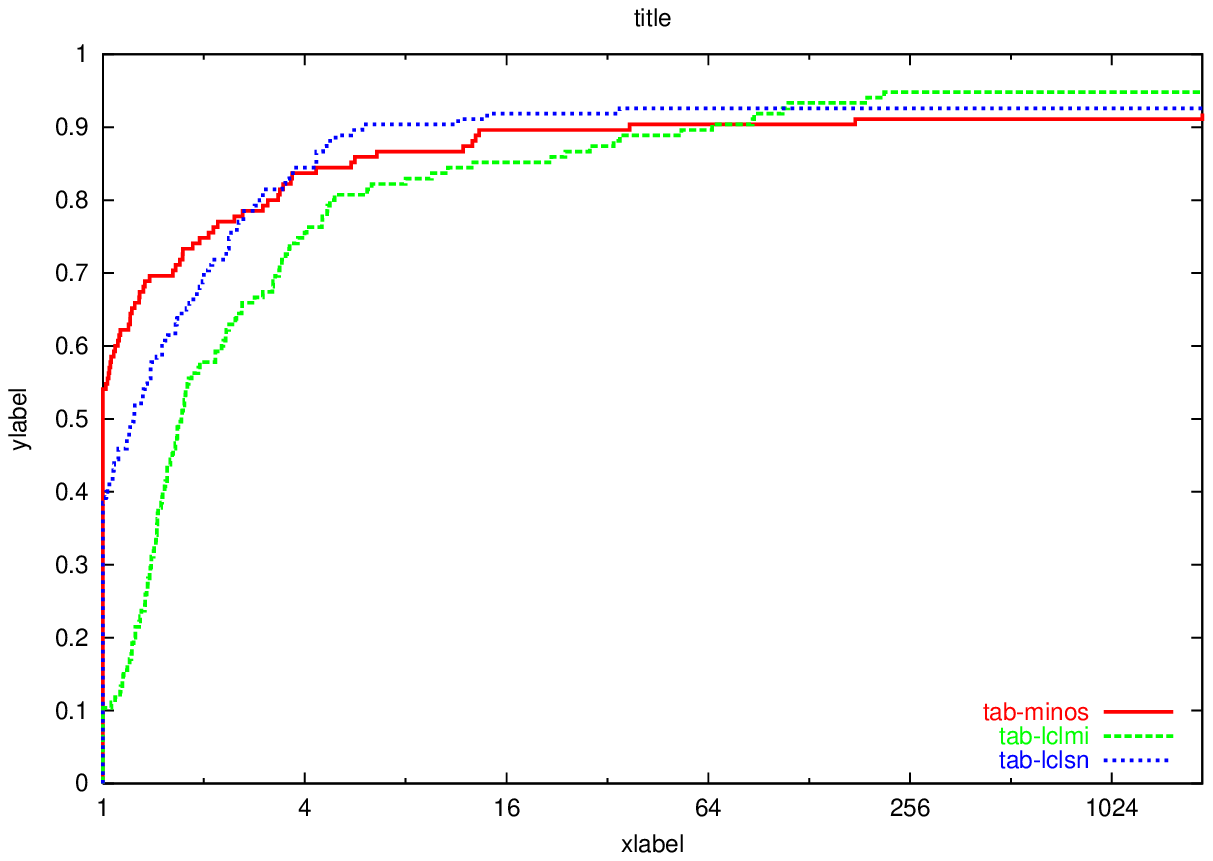}}\\
  \subfigure[Minor iterations]%
  {\includegraphics[height=2.3in]{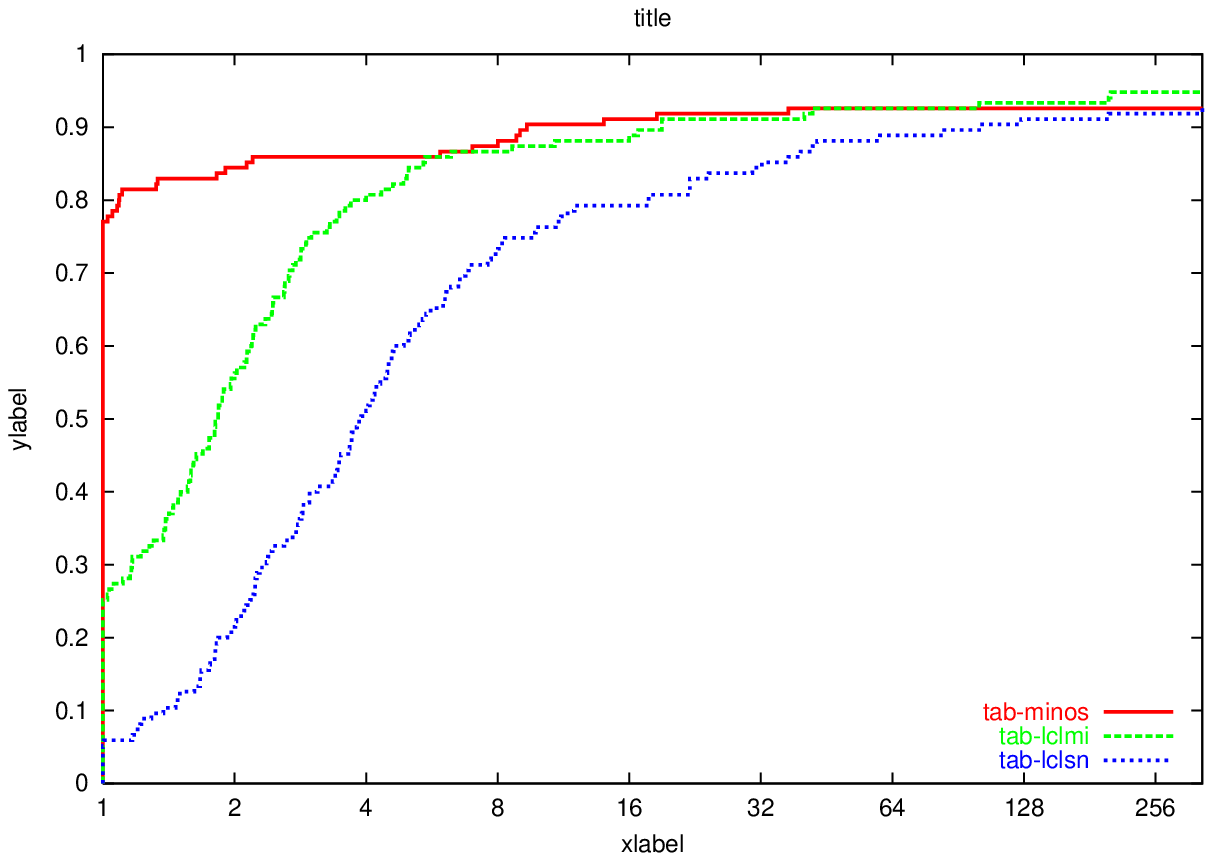}}\\
  \subfigure[Major iterations]%
  {\includegraphics[height=2.3in]{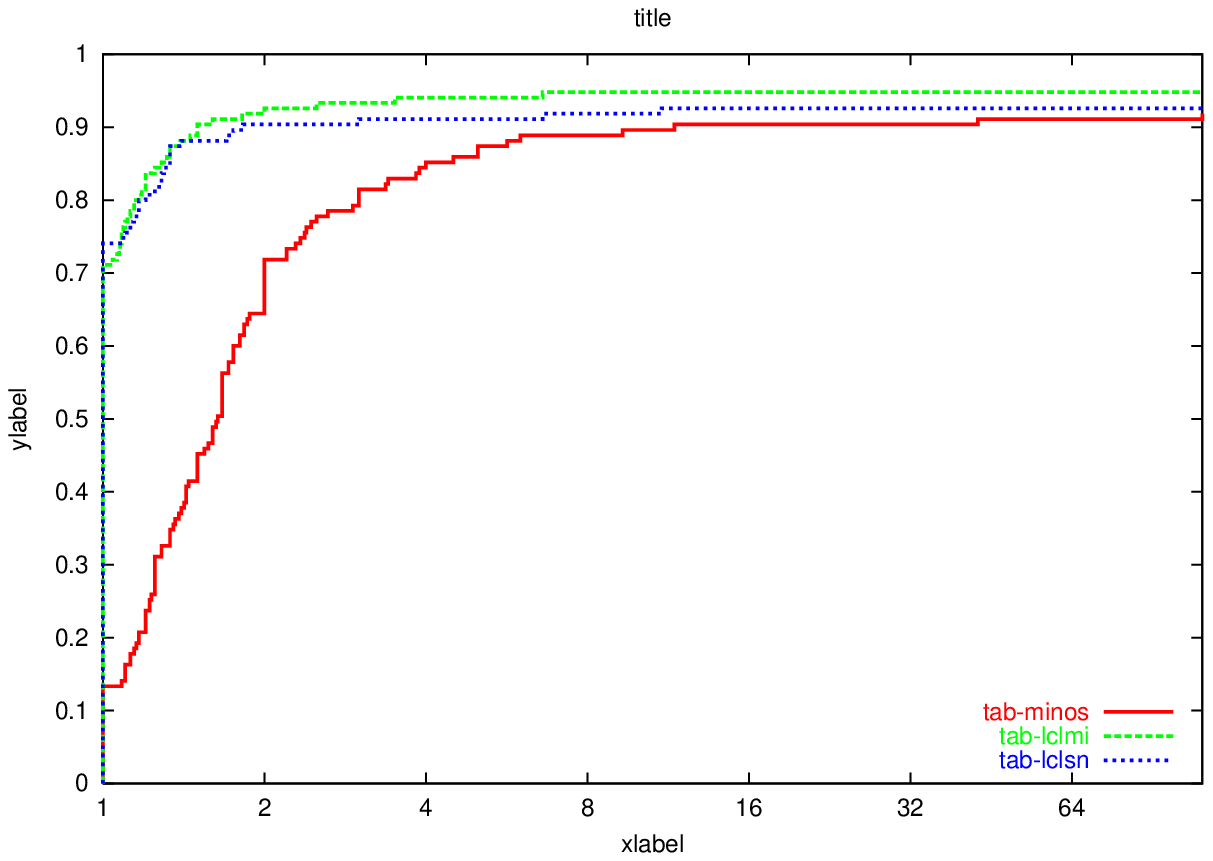}}
  \caption{Performance profiles.   The vertical axes represent the
    percentage of problems successfully solved within a factor $\tau$
    of the best solver.  The horizontal axes are based on a log scale.
    Performance profiles are shown for the number of nonlinear
    function evaluations, minor iterations, and major iterations.  The
    profiles include the results of the 135 selected test problems.}
\end{figure}

Figure~\ref{fig:perfprof} shows the \emph{performance profiles}, as
described by Dolan and Mor\'{e}~\cite{DM01}, of the two versions of
\LCLOPT\ (the dotted and dashed lines) and \MINOS\ (the solid line).
The three charts of that figure show performance profiles for the
total number of nonlinear function evaluations, minor iterations, and
major iterations.  All the problems selected from the \COPS,
Hock-Schittkowski, and \CUTE\ test suites are included in each
profile.  The performance profiles describe the percentage of problems
successfully solved (the vertical axes) within a factor $\tau$ of the
best-performing solver (the horizonal axes).

By all measures, \LCLOPT, using \MINOS\ to solve the subproblems,
successfully solves the largest proportion of problems and proves to
be the most reliable method.  Compared with \MINOS, \LCLOPT\ tends to
require more minor iterations (a measure of total computational work)
but fewer major iterations to reach a solution.  We comment further on
this fact in \Sec\ref{sec:import-early-term}.

\renewcommand{\subfigcapskip}{0pt} \renewcommand{\subfigtopskip}{0pt}

For all problems that can vary in the number of constraints and
variables, we describe their dimensions.  The following heads are used
in Tables~\ref{tab:COPSsizes} and~\ref{tab:CUTEsizes}.

\begin{center}
\small
\begin{tabular}{@{}ll@{}}
   \toprule
   Head   & Dimension
\\ \midrule
   $m$    & Constraints (linear and nonlinear)
\\ $m_c$  & Nonlinear constraints
\\ $n$    & Variables
\\ $n_c$  & Variables appearing nonlinearly in $c$
\\ $n_f$  & Variables appearing nonlinearly in $f$
\\ \bottomrule
\end{tabular}
\end{center}

\subsection{Default parameters}                  \label{sec:num-params}

Figure~\ref{fig:LCspecs} shows the options files that \LCLOPT\ uses
for the LC solvers.  These are fixed for all subproblems.  Separately,
at each major iteration, \LCLOPT\ sets the parameter \opt{Optimality
  Tolerance} in \MINOS\, and the parameter \opt{Major Optimality
  Tolerance} in \SNOPT.  These are equivalent to the subproblem
optimality tolerance $\omega\k$ (\Sec\ref{sec:implem-early-term}
outlines the method for choosing this parameter).

Each test problem supplies a default starting point.  This point is
used as $\xtilde$ in the proximal-point problem (see
\Sec\ref{sec:implem-init}).  The initial vector of multiplier
estimates $y_0$ is set to zero.

Both \MINOS\ and \SNOPT\ provide the option to reuse a quasi-Newton
approximation of a Hessian from a previous solve: \MINOS\ approximates
the reduced Hessian; \SNOPT\ approximates the full Hessian.  We take
advantage of this feature for all iterations $k=2,3,4,\ldots$ by
setting the \MINOS\ and \SNOPT\ options \opt{Start = `Hot'}.

The parameters used by Algorithm~\ref{algo:sLCL} are set as follows.
The upper and lower bounds of the elastic penalty parameters are
$\sigmab = 1$ and $\sigmaB = 10^4$.  The initial elastic weight
is $\sigma_0=10^2$. (Normally, \LCLOPT\ scales this quantity by
$1+\infnorm{y_0}$, but the scaling has no effect for these test runs
because $y_0\equiv0$.)  The penalty scaling factors are $\tau_{\rho} =
100^{0.5}$ and $\tau_{\sigma}=10$.  As suggested in~\cite{CGT91b}, we
set $\alpha = 0.1$ and $\beta=0.9$.  The initial penalty parameter is
$\rho_0=10^{5/2}/m_c$, where $m_c$ is the number of nonlinear
constraints.  The final optimality and feasibility tolerances are
$\omega_*=\eta_*=10^{-6}$.  The initial optimality and feasibility
tolerances are $\omega_0=10^{-3}$ ($=\sqrt{\omega_*}$) and $\eta_0=1$.

In all cases, default options, with the exception of \opt{Major
  Iterations 500} and \opt{Superbasics Limit 2000}, are used for the
\MINOS\ benchmarks.

\begin{figure}[ht]
\small
\begin{minipage}{\textwidth}
\centering
\begin{boxedminipage}[b]{.45\textwidth}
\begin{verbatim}

 BEGIN LCL SUBPROBLEM
    Scale option                0
    Superbasics limit        2000
    Iterations               5000
    Feasibility tol        1.0e-6
 END LCL SUBPROBLEM
\end{verbatim}
\centering
(a) The \MINOS\ specs file
\end{boxedminipage}
\hfill
\begin{boxedminipage}[b]{.45\textwidth}
\begin{verbatim}

 BEGIN LCL SUBPROBLEM
    Scale option                0
    Superbasics limit        2000
    Iterations               5000
    Major iterations         1000
    Minor iterations          500
    Minor feasibility tol  1.0e-6
    Minor optimality  tol  2.5e-7
 END LCL SUBPROBLEM
\end{verbatim}
\centering
(b) The \SNOPT\ specs file
\end{boxedminipage}
\end{minipage}
  \caption[The fixed optional parameters for every subproblem solve]
  {The fixed optional parameters for every subproblem solve.  The
    optimality tolerance $\omega\k$ is specified by \LCLOPT\ for each
    $k$.}
  \label{fig:LCspecs}
\end{figure}

\subsection{The COPS test problems}                \label{sec:num-cops}

The \COPS~2.0 collection~\cite{DM00} comprises 17 problems.  Five
problems are excluded for the following reasons:
\begin{itemize}
\item 3 problems are unconstrained: \prob{bearing}, \prob{minsurf},
  and \prob{torsion};
\item 2 problems cause system errors when called using the \AMPL\ MEX
  interface: \prob{glider} and \prob{marine}.
\end{itemize}
The dimensions of the \COPS\ test problems can be adjusted.  In all
cases, the solvers were applied to the largest version of the problem
(as specified by the \AMPL\ model) that would not cause the system to
age memory to disk.  Table~\ref{tab:COPSsizes} summarizes the
dimensions of the selected problems.

\begin{table}[htb]
\caption{Dimensions: The 12 selected \COPS\ test problems}
\label{tab:COPSsizes}
\centering
\small
\begin{tabular}{@{}lrrrrr@{}}
  \toprule
  Problem & $m$ & $m_c$ &$n$    & $n_c$ & $n_f$\\
  \midrule
 camshape   &  1604 &   801 &   800 &   800 &     0\\
 catmix     &  1603 &  1600 &  2403 &  2403 &     0\\
 chain      &   204 &     1 &   402 &   201 &   402\\
 channel    &   800 &   400 &   800 &   800 &     0\\
 elec       &   201 &   200 &   600 &   600 &   600\\
 gasoil400  &  4004 &  3200 &  4003 &  4003 &   202\\
 marine     &  1208 &   800 &  1215 &  1215 &   344\\
 methanol   &  2406 &  1800 &  2405 &  1605 &  1670\\
 pinene     &  4006 &  3000 &  4005 &  2405 &  2469\\
 polygon    &  1377 &  1225 &   100 &   100 &   100\\
 robot      &  2414 &  2400 &  3611 &  3209 &     0\\
 rocket     &  2409 &  1200 &  1605 &  1605 &     0\\
 steering   &  2011 &  1600 &  2007 &  1204 &     0\\
  \bottomrule
\end{tabular}
\end{table}

As shown in Table~\ref{tab:COPSsumm},
the version of \LCLOPT\ using \MINOS\ for the subproblems solved all
12 problems to first-order optimality.  The version using \SNOPT\ 
solved 11 problems to first-order optimality; the exception was
\prob{robot}, which it declared as having infeasible nonlinear
constraints.  \MINOS\ solved 10 of the 12 problems to optimality; it
declared \prob{steering} an infeasible problem, and it terminated the
solution of \prob{elec} because of excessive iterations.  Feasible
points exist for all of the test problems chosen, so we consider all
declarations of infeasibility to be errors.

\begin{table}[htb]
  \centering
  \small
  \caption{Summary: The 12 selected \COPS\ test problems}
  \label{tab:COPSsumm}
  \begin{tabular}{@{}ld{2}d{2}d{2}@{}}
   \toprule
         & \mcol{2}{c}{\LCLOPT}&
\\ \cmidrule(r){2-3}
         & \mcol{1}{@{}c}{(\MINOS)}
         & \mcol{1}{@{}c}{(\SNOPT)}
         & \mcol{1}{c@{}}{\MINOS}
\\ \midrule
   Optimal                           & 12 & 11 & 10
\\ False Infeasibility               &    &  1 &  1
\\ Terminated by iteration limit     &    &    &  1
\\ \midrule
   Major iterations                  &   118 &    179 &   380 
\\ Minor iterations                  & 53950 & 147518 & 61388 
\\ Function evaluations              & 53081 &  11014 & 63701
\\ \bottomrule
  \end{tabular}
\end{table}

We note that different local optima appear to have been found for
problems \prob{camshape, methanol, polygon}, and \prob{rocket}.  An
excessive number of minor iterations were required by \LCLOPT\ on
\prob{catmix, elec}, and \prob{robot} with \SNOPT\ as its subproblem
solver.  Especially during early major iterations, \SNOPT\ was unable
to solve the LC subproblems to the required optimality tolerance
within the 5000 iteration limit.  Rather than terminate with an error
message, \LCLOPT\ forces \SNOPT\ to keep working on the same
subproblem until it returns a solution within the required optimality
tolerance.  In practice, a different strategy would be adopted, but
our goal here is to test the robustness of the outer iterations (the
stabilized LCL method), not the robustness of the subproblem solvers.

\subsection{The Hock-Schittkowski test problems}     \label{sec:num-hs}

The HS test suite contains 86 nonlinearly constrained
problems~\cite{HS81}.  These are generally small and dense problems.
We exclude 5 problems from this set for the following reasons:
\begin{itemize}
\item 3 problems are not smooth: \prob{hs67, hs85}, and \prob{hs87};
\item 2 problems require external functions: \prob{hs68} and
  \prob{hs69}.
\end{itemize}

Both versions of \LCLOPT\ solved the same 80 problems to first-order
optimality, but both declared \prob{hs109} infeasible.  \MINOS\ solved
80 problems to first-order optimality but declared \prob{hs93}
infeasible.

\begin{table}[htb]
  \small
  \centering
  \caption{Summary: The 81 selected Hock-Schittkowski test problems}
  \label{tab:HSsumm}
  \begin{tabular}{@{}ld{2}d{2}d{2}@{}}
   \toprule
         & \mcol{2}{c}{\LCLOPT}&
\\ \cmidrule(r){2-3}
         & \mcol{1}{@{}c}{(\MINOS)}
         & \mcol{1}{@{}c}{(\SNOPT)}
         & \mcol{1}{c@{}}{\MINOS}
\\ \midrule
   Optimal                           & 80    & 80     & 80
\\ False infeasibility               &  1    &  1     &  1
\\ \midrule
   Major iterations                  &   654 &    648 &  1160 
\\ Minor iterations                  &  7415 &  25290 & 10111 
\\ Function evaluations              & 12269 &  14712 & 27127
\\ \bottomrule
  \end{tabular}
\end{table}

On \prob{hs13}, all the solvers reached different solutions.  However,
the linear independence constraint qualification does not hold at the
solution of this problem---this violates the required assumptions for
both \LCLOPT\ and \MINOS.

Recall that \LCLOPT\ and \MINOS\ use only first derivatives and hence
may not necessarily converge to local solutions of a problem.  For
example, \LCLOPT\ (in both versions) converged to a known local
solution of \prob{hs16}, but \MINOS\ converged to some other
first-order point.  In contrast, \MINOS\ converged to the known local
solutions of \prob{hs97} and \prob{hs98}, while \LCLOPT\ (in both
versions) converged to other first-order points.  Similar differences
exist for problems \prob{hs47} and \prob{hs77}.

\subsection{A selection of CUTE test problems}     \label{sec:num-cute}

With the \texttt{select} utility~\cite{BCGT95}, we extracted from the
\CUTE\ test suite dated September 7, 2000, problems with the following
characteristics (\verb+*+ is a wild-card character):

\begin{center}
\begin{boxedminipage}{\textwidth}
\small
\begin{verbatim}

 Objective function type          : *
 Constraint type                  : Q O (quadratic, general nonlinear)
 Regularity                       : R   (smooth)
 Degree of available derivatives  : 1   (first derivatives, at least)
 Problem interest                 : M R (modeling, real applications)
 Explicit internal variables      : *
 Number of variables              : *
 Number of constraints            : *
\end{verbatim}
\vspace*{0.5\baselineskip}
\end{boxedminipage}
\end{center}

\null
\noindent
These criteria yield 108 problems.  We exclude 66 problems from this
set for the following reasons:
\begin{small}
\begin{itemize}
\item 33 problems do not have \AMPL\ versions: {\it 
car2,
c-reload, 
dembo7, 
drugdis, 
durgdise,
errinbar, 
junkturn, 
leaknet,
lubrif, 
mribasis, 
nystrom5, 
orbit2, 
reading4, 
reading5, 
reading6, 
reading7, 
reading8, 
reading9, 
rotodisc,
saromm, 
saro, 
tenbars1, 
tenbars2, 
tenbars3, 
tenbars4, 
trigger,
truspyr1, 
truspyr2, 
zamb2, 
zamb2-8, 
zamb2-9, 
zamb2-10}, and 
{\it zamb2-11};
\item 21 problems cause system errors when evaluated either by the \AMPL\ 
  MEX interface or by \MINOS\ (when invoked from \AMPL): {\it
brainpc2,
brainpc3,
brainpc4, 
brainpc5, 
brainpc6, 
brainpc7, 
brainpc8, 
brainpc9,
bratu2dt, 
cresc132, 
csfi1, 
csfi2, 
drcav1lq, 
drcav2lq, 
drcav3lq,
kissing, 
lakes, 
porous1, 
porous2, 
trainf}, and
{\it trainh};
\item The \AMPL\ versions of 12 problems are formulated with no
  nonlinear constraints: {\it 
drcavty1, 
drcavty2, 
drcavty3, 
flosp2hh,
flosp2hl, 
flosp2hm, 
flosp2th, 
flosp2tl,\\
flosp2tm, 
methanb8,
methanl8}, and
{\it res}.
\end{itemize}
\end{small}
The dimensions of 17 of the remaining 42 problems can be adjusted.  In
all cases, the solvers were applied to the largest problem versions
that would not cause the system to page memory to disk.
Table~\ref{tab:CUTEsizes} summarizes the dimensions of the selected
problems that can vary in size.

\begin{table}[htb]
\small
\caption{Dimensions of the variable-size \CUTE\ test problems}
\label{tab:CUTEsizes}
\centering
\begin{tabular}{@{}lrrrrr@{}}
  \toprule
  Problem & $m$ & $m_c$ &$n$    & $n_c$ & $n_f$\\
  \midrule
 bdvalue &  1000 &  1000 &  1000 &  1000 &     0\\
 bratu2d &  4900 &  4900 &  4900 &  4900 &     0\\
 bratu3d &   512 &   512 &   512 &   512 &     0\\
 cbratu2d&   882 &   882 &   882 &   882 &     0\\
 cbratu3d&  1024 &  1024 &  1024 &  1024 &     0\\
 chandheq&   100 &   100 &   100 &   100 &     0\\
 chemrcta&  2000 &  1996 &  2000 &  1996 &     0\\
 chemrctb&  1000 &   998 &  1000 &   998 &     0\\
 clnlbeam&  1001 &   500 &  1499 &   499 &  1000\\
 hadamard&   257 &   128 &    65 &    64 &    65\\
 manne   &   731 &   364 &  1094 &   364 &   729\\
 reading1&  5001 &  5000 & 10001 & 10000 & 10000\\
 reading3&   103 &   101 &   202 &   202 &   202\\
 sreadin3&  5001 &  5000 & 10000 &  9998 &  9998\\
 ssnlbeam&    21 &    10 &    31 &    11 &    22\\
 svanberg&  1001 &  1000 &  1000 &  1000 &  1000\\
 ubh5    & 14001 &  2000 & 19997 &  6003 &     0\\
 \bottomrule
\end{tabular}
\end{table}

The version of \LCLOPT\ using \MINOS\ solved 36 of 42 problems to
first-order optimality, while the version using \SNOPT\ solved 34
problems to first-order optimality.  \MINOS\ solved 34 problems to
first-order optimality.  Table~\ref{tab:CUTEsumm} summarizes these
results.  We note that \LCLOPT, in one of its two versions, solves
every problem except \prob{heart6}, which it declares infeasible.
With the exception of \prob{cresc50}, \LCLOPT\ with \SNOPT\ does not
seem to suffer (on successful solves) from excessive minor iterations
resulting from subproblem restarts, as it does on the \COPS\ problems.

\begin{table}[htb]
  \small
  \centering
  \caption{Summary: The 42 selected \CUTE\ test problems}
  \label{tab:CUTEsumm}
  \begin{tabular}{@{}ld{2}d{2}d{2}@{}}
   \toprule
         & \mcol{2}{c}{\LCLOPT}&
\\ \cmidrule(r){2-3}
         & \mcol{1}{@{}c}{(\MINOS)}
         & \mcol{1}{@{}c}{(\SNOPT)}
         & \mcol{1}{c@{}}{\MINOS}
\\ \midrule
   Optimal                           & 36     & 34     &    34
\\ False infeasibility               &  4     &  3     &     3
\\ Terminated by iteration limit     &  1     &  3     &     1
\\ Terminated by superbasics limit   &        &  1     & 
\\ Unbounded/badly scaled            &        &        &     3
\\ Final point cannot be improved    &  1     &  1     &     1
\\ \midrule
   Major iterations                  &    400 &    368 &  1149 
\\ Minor iterations                  &  70476 & 278162 & 29021 
\\ Function evaluations              &  59216 &  57732 & 53069
\\ \bottomrule
  \end{tabular}
\end{table}


\section{Conclusions}                                   \label{sec:con}

The stabilized LCL method developed in this paper is a generalization
of the augmented Lagrangian methods discussed in
\Sec\ref{sec:sLCL-algo} and it shares the strengths of its
predecessors: it is globally convergent (the BCL advantage) and it has
fast local convergence (the LCL advantage).  The $\ell_1$-penalty
function brings the two together.  Because the stabilized LCL method
operates in a reduced space given by the linearized constraints (like
the LCL method), it does not suffer from the ill-conditioning effects
that can plague BCL methods.

\subsection{Importance of early termination}\label{sec:import-early-term}

The numerical results presented in \Sec\ref{sec:num} demonstrate that
\MINOS\ successfully solved many of the test problems using relatively
few minor iterations.  \MINOS\ terminates its progress on each of its
subproblems after 40 iterations (to avoid a refactorization of the
current basis, which by default occurs every 50 iterations).  In
contrast, \LCLOPT\ attempts to constrain the subproblem iterations by
means of an initially loose optimality tolerance (we set
$\omega_0=\sqrt{\omega_*}$ for the runs shown in \Sec\ref{sec:num}).
A potential weakness of this approach vis \`{a} vis \MINOS\ is that
there is no a priori bound on the number of subproblem iterations.
\MINOS's aggressive (and heuristic) strategy seems effective in
keeping the total minor iteration counts low.  This property is
particularly important during the early major iterations, when the
current solution estimates are poor.

It may be possible to emulate the \MINOS\ strategy and still satisfy
the requirement that the subproblem optimality tolerances
$\omega\k$ converge to zero (\cf~Lemma~\ref{le:penalty}).  For
example, \LCLOPT\ might truncate the subproblem solutions after a
fixed number of iterations, and only gradually increase the iteration
limit on successive major iterations.  Especially during early major
iterations, such a strategy may keep the accumulated number of
subproblem iterations small.  During later major iterations, the
strategy would still ensure that the subproblem solver returns
solutions within the prescribed tolerance $\omega\k$.

On the other end of the performance spectrum lies the issue of
recovering LCL's fast local convergence rate under inexact solves
(\cf~\Sec\ref{sec:sLCL-rates}).  Br\"{a}uninger~\cite{Bra81}\ proves
that the quadratic convergence rate of Robinson's method is retained
when $\omega\k$ is reduced at a rate $O(\norm{F(x\k,y\k,z\k)}^2)$
(\cf~Theorem~\ref{th:lclrates}).  The first-order KKT
conditions~\eqref{eq:58} for the LCL subproblem can be expressed as
\begin{equation}
  \label{eq:26}
  \begin{pmatrix}
    \Hess{xx}\AugL\k(x\k)  &  J\k\T
\\  J\k                    &
  \end{pmatrix}
  \begin{pmatrix}
    p \\ -y
  \end{pmatrix}
  + O(\norm{p}^2) =
  \begin{pmatrix}
    -g\k + J\k\T y\k
\\  -c\k
  \end{pmatrix},
\end{equation}
where $p\defd x-x\k$, and a first-order Taylor expansion was used to
derive the residual term $O(\norm{p}^2)$. (We have ignored bound
constraints for the moment.  Robinson~\cite{Rob72,Rob74}\ shows that
the correct active set is identified by the subproblems near a
solution.)  The nonlinear equations~\eqref{eq:26} are closely related
to the linear equations that would be derived from applying Newton's
method to~\eqref{eq:58} (again, ignoring bound constraints).  In that
case, the theory from inexact Newton methods (Dembo
\etal~\cite{DES82}) predicts that the quadratic convergence rate is
recovered when the residual error is reduced at the rate
$O(\norm{F(x\k,y\k,z\k)})$.  The similarity between~\eqref{eq:26} and
the Newton equations hints at the possibility of recovering the
quadratic convergence rate of the LCL and stabilized LCL methods by
reducing $\omega\k$ at the rate $O(\norm{F(x\k,y\k,z\k)})$.  We note,
however, that stronger assumptions may be needed on the smoothness of
the nonlinear functions.  This issue deserves more study.

\subsection{Keeping the penalty parameter small}

Preliminary experimentation reveals that a small penalty parameter
$\rho\k$ can significantly reduce the difficulty of each subproblem
solve.  BCL methods require that $\rho\k$ be larger than some
threshold value $\rhobar$.  In contrast, LCL methods can converge when
$\rho\k\equiv0$ if they are started near a solution (see
\Sec\ref{th:lclrates}).

The challenge here is to find a strategy that can keep $\rho\k$ small
or reduce it without destabilizing the method.  A tentative strategy
might be to reduce $\rho\k$ only finitely many times.  This approach
does not violate the hypotheses of Lemma~\ref{le:penalty}, and may be
effective in practice.  A form of this strategy was used for the runs
shown in \Sec\ref{sec:num}.

\subsection{A second-derivative LC solver}        \label{sec:2nd-order}

We prove in \Sec\ref{sec:sLCL-second} that the stabilized LCL method
will converge to second-order stationary points if the LC subproblems
are solved to second-order points (for example, by using a
second-derivative LC solver).  In practice, however, a
second-derivative LC solver may be most useful as a means of reducing
the overall computational work required by the stabilized LCL method.

The stabilized LCL method is largely independent of the method in
which its subproblems are solved.  An LC solver using second
derivatives is likely to require fewer iterations (and hence less
computational work) for the solution of each of the subproblem.  We
would expect the number of required major iterations to remain
constant if each subproblem solution is computed to within the
prescribed tolerance $\omega\k$.  However, we would expect to
\emph{reduce} the number of required major iterations if a \MINOS-like
strategy is used to terminate the subproblems (see
\Sec\ref{sec:import-early-term}).  Over the same number of iterations,
a subproblem solver using second derivatives may make more progress
toward a solution than a first-derivative solver.

Any future implementation of the stabilized LCL method would ideally
be flexible enough to allow for a variety of solvers to be used for
the LC subproblems.  The choice of the subproblem solver could then be
guided by the characteristics of the optimization problem at hand.  In
particular, the advent of automatic differentiation makes second
derivatives increasingly available for certain problem classes, e.g.,
within recent versions of GAMS and AMPL, and for more general
functions defined by Fortran or C code, notably ADIFOR and ADIC
(Bischof \etal~\cite{BR97,BCH98}).  These may be used by SQP and
interior methods for nonlinearly constrained (NC) problems (e.g., LOQO
Vanderbei~\cite{SV99}).  Certain theoretical challenges might be
avoided, however, by developing specialized second-derivative LC
solvers.  Such LC solvers could be extended readily to general NC
problems by incorporating them into the stabilized LCL algorithm.

\bibliography{lcl_ref}

\vfill
\begin{flushright}
\scriptsize
\framebox{\parbox{2.4in}{The submitted manuscript has been created
by the University of Chicago as Operator of Argonne
National Laboratory ("Argonne") under Contract No.\
W-31-109-ENG-38 with the U.S. Department of Energy.
The U.S. Government retains for itself, and others
acting on its behalf, a paid-up, nonexclusive, irrevocable
worldwide license in said article to reproduce,
prepare derivative works, distribute copies to the
public, and perform publicly and display publicly, by or on
behalf of the Government.}}
\normalsize
\end{flushright}

\end{document}